\documentclass{gtart_h}  

\def\ifplaintex{\expandafter\ifx\csname documentclass\endcsname\relax}

\def\ifplaintex{\expandafter\ifx\csname documentclass\endcsname\relax}

\ifplaintex 
\else
\fi

\expandafter\ifx\csname epsfbox\endcsname\relax\input epsf\fi

\def\gt{{\mathsurround=0pt\it $\cal G\mskip-2mu$eometry \&\ 
$\cal T\!\!$opology}}        

\def\gtp{{\mathsurround=0pt\it $\cal G\mskip-2mu$eometry \&\ 
$\cal T\!\!$opology $\cal P\!$ublications}}  

\def\lognumber#1{\def\thelognumber{#1}}
\def\volumenumber#1{\def\thevolumenumber{#1}}
\def\papernumber#1{\def\thepapernumber{#1}}
\def\volumeyear#1{\def\thevolumeyear{#1}}

\def\pagenumbers#1#2{\def\startpage{#1}\def\finishpage{#2}}
\def\published#1{\def\publishdate{#1}}
\def\proposed#1{\def\theproposer{#1}}
\def\seconded#1{\def\theseconders{#1}}
\def\received#1{\def\receiveddate{#1}}
\def\revised#1{\def\reviseddate{#1}}
\def\accepted#1{\def\accepteddate{#1}}

\def\asciiemail#1{\def\theasciiemail{#1}}

\long\def\asciiabstract#1{\long\def\theasciiabstract{#1}}

\def\shorttitle#1{\def\theshorttitle{#1}}

\let\\\par\let\thelognumber\relax
\let\thevolumenumber\relax\let\thepapernumber\relax
\let\thevolumeyear\relax\let\thesamplenumber\relax\let\startpage\relax
\let\finishpage\relax\let\publishdate\relax\let\receiveddate\relax
\let\reviseddate\relax\let\accepteddate\relax\let\theasciititle\relax
\let\theasciiauthors\relax
\let\theasciiabstract\relax
\let\theasciiemail\relax\let\theshortauthors\relax\let\theshorttitle\relax

\long\def\maketitlep{   

\count0=\startpage

\gt\hfill      
\hbox to 77pt{\vbox to 0pt{\vglue -15pt\epsfbox{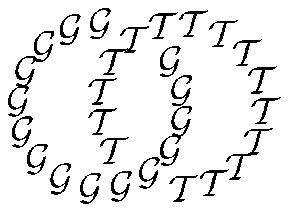}\vss}\hss}
\break
{\small\ifx\thesamplenumber\relax 
Volume \else Sample
\fi\thevolumenumber\ (\thevolumeyear)
\startpage--\finishpage\nl
Published: \publishdate}
\vglue 0.5truein plus 0.4fil minus 0.1truein

{\parskip=0pt\leftskip 0pt plus 1fil\def\\{\par\smallskip}{\ifplaintex\large
\else\Large\fi\bf\thetitle}\par\medskip}   

\vglue 0pt plus 0.1fil 

{\parskip=0pt\leftskip 0pt plus 1fil\def\\{\par}{\sc\theauthors}
\par\medskip}

\vglue 0pt plus 0.1fil 

{\small\parskip=0pt\let\newline\\
{\leftskip 0pt plus 1fil\def\\{\par}{\sl\theaddress}\par}
\expandafter\ifx\theemail\relax    
\relax\else\vglue 5pt plus 0.02fil minus 2pt\def\\{\stdspace{\rm 
and}\stdspace} 
\cl{Email:\stdspace\tt\theemail}\fi
\ifx\theurl\relax                  
\relax\else\vglue 5pt plus 0.02fil minus 2pt\def\\{\stdspace{\rm 
and}\stdspace}
\cl{URL:\stdspace\tt\theurl}\fi\par}

\vglue 7pt plus 0.3fil minus 3pt

{\bf Abstract}
\vglue 5pt plus 0.1fil minus 2pt

\theabstract

\vglue 7pt plus 0.3fil minus 3pt

{\bf AMS Classification numbers}\quad Primary:\quad \theprimaryclass

Secondary:\quad \thesecondaryclass

\vglue 5pt plus 0.3fil minus 2pt

{\bf Keywords:}\quad \thekeywords

\vglue 10pt plus 0.5fil minus 5pt

{\small  Proposed: \theproposer\hfill Received: \receiveddate\nl
Seconded: \theseconders\hfill 
\ifx\reviseddate\relax                         
Accepted: \accepteddate                        
\else
Revised: \reviseddate                          
\fi}
\eject
}       

\font\phead=cmsl9 scaled 950
\font=cmsl9 scaled 1050
\font\pnum=cmbx10 scaled 913
\font\lnum=cmbx10 
\font\pfoot=cmsl9 scaled 950
\font=cmsl9 scaled 1050
\ifplaintex
\headline{\vbox to 0pt{\vskip -4.5mm\line{\small\phead\ifnum
\count0=\startpage ISSN 1364-0380 (on line)
1465-3060 (printed) \hfill {\pnum\folio}\else\ifodd\count0\def\\{ }%
\ifx\theshorttitle\relax\thetitle\else\theshorttitle\fi\hfill{\pnum\folio}
\else\def\\{ and }{\pnum\folio}\hfill\ifx\theshortauthors\relax\theauthors
\else\theshortauthors\fi\fi\fi}\vss}}
\footline{\vbox to 0pt{\vglue 0mm\line{\small\pfoot\ifnum\count0=\startpage
\copyright\ \gtp\hfill\else
\gt, Volume \thevolumenumber\ (\thevolumeyear)\hfill\fi}\vss
}}
\else
\makeatletter
\def\@oddhead{{\small\lhead\ifnum\count0=\startpage ISSN 1364-0380 (on line)
1465-3060 (printed) \hfill {\lnum\number\count0}\else\ifodd\count0
\def\\{ }\ifx\theshorttitle\relax \thetitle \else\theshorttitle\fi\hfill
{\lnum\number\count0}\else\def\\{ and }{\lnum\number\count0}
\hfill\ifx\theshortauthors\relax 
\theauthors\else\theshortauthors\fi\fi\fi}}\def\@evenhead{\@oddhead}
\def\@oddfoot{\small\lfoot\ifnum\count0=\startpage\copyright\ \gtp\hfill\else
\gt, Volume \thevolumenumber\ (\thevolumeyear)\hfill\fi}
\def\@evenfoot{\@oddfoot}
\makeatother
\fi

\newwrite\gtoutfile
\long\gdef\makeheadfile{  
{\def\\{, }\def\s{ }
\immediate\openout\gtoutfile head.xxx
\immediate\write\gtoutfile{Proxy-for: \ifx\theasciiauthors\relax
\theauthors\else\theasciiauthors\fi\s<\ifx\theasciiemail\relax\theemail\else\theasciiemail\fi>}
\immediate\write\gtoutfile{\noexpand\\}
\immediate\write\gtoutfile{Authors: \ifx\theasciiauthors\relax
\theauthors\else\theasciiauthors\fi}
{\def\\{ }\immediate\write\gtoutfile{Title: \ifx\theasciititle\relax
\thetitle\else\theasciititle\fi}}
\immediate\write\gtoutfile{Subj-class: GT or SG or MG etc}
\immediate\write\gtoutfile{MSC-class: \theprimaryclass\ifx\thesecondaryclass\relax\else, \thesecondaryclass\fi}
\immediate\write\gtoutfile{Journal-ref: Geom. Topol. \thevolumenumber
(\thevolumeyear) \startpage-\finishpage}
\immediate\write\gtoutfile{Comments: Published by Geometry and Topology at}
\immediate\write\gtoutfile{\s\s http://www.maths.warwick.ac.uk/gt/GTVol\thevolumenumber/paper\thepapernumber.abs.html}
\immediate\write\gtoutfile{\noexpand\\}
\immediate\write\gtoutfile{}
\ifx\theasciiabstract\relax
\immediate\write\gtoutfile{\theabstract}\else
\immediate\write\gtoutfile{\theasciiabstract}\fi
\immediate\write\gtoutfile{}
\immediate\write\gtoutfile{\noexpand\\}
\immediate\write\gtoutfile{}
\immediate\closeout\gtoutfile}}  

\def\maketitlepage{\maketitlep\makeheadfile}
\let\maketitle\maketitlepage

\lognumber{317}
\volumenumber{8}\papernumber{15}\volumeyear{2004}
\pagenumbers{611}{644}
\received{10 April 2003}
\revised{8 April 2004}
\published{12 April 2004}
\accepted{19 December 2004}
\proposed{Benson Farb}
\seconded{Martin Bridson, Steven Ferry}

\usepackage{xypic, amssymb, epic, buxmath}

\newif\ifXYpic
  \newdir{ >}{%
    !/10pt/@{}*@{>}
  }
  \XYpictrue
\renewcommand{\monorightarrow}{\rightarrowtail}

\DeclareSymbolFont{largesymsup}{U}{bmexs}{m}{n}
\DeclareMathSymbol{\bigast}{\mathop}{largesymsup}{'100}
\DeclareMathSymbol{\bigudot}{\mathop}{largesymsup}{'102}
\DeclareMathSymbol{\bigplus}{\mathop}{largesymsup}{'110}
\DeclareMathSymbol{\bigtimes}{\mathop}{largesymsup}{'112}

\newenvironment{refferences}{%
  \newcommand{\au}[1]{\textbf{##1},}%
  \newcommand{\ed}[1]{\textbf{##1},}%
  \newcommand{\ti}[1]{\emph{##1},}%
  \newcommand{\edition}[1]{##1 }%
  \newcommand{\lo}[2]{%
    ##1%
    \begingroup
      \def\help{##2}%
      \def\empty{}%
      \ifx\empty\help\else
        ~\help
      \fi
    \endgroup
  }%
  \newcommand{\book}[2]{%
    \bibitem{##2}%
  }%
  \newcommand{\report}[2]{%
    \bibitem{##2}%
  }%
  \newcommand{\article}[2]{%
    \bibitem{##2}%
  }%
  \newcommand{\thesis}[2]{%
    \bibitem{##2}%
  }%
  \newcommand{\notes}[2]{%
    \bibitem{##2}%
  }%
  \newcommand{\preprint}[2]{%
    \bibitem{##2}%
  }%
  \thebibliography
}{%
  \endthebibliography
}

\renewcommand{\notion}{\emph}
  \newcommand{\xnotion}[1]{#1}
  \newenvironment{sketch}{%
    \begin{proof}[Sketch of proof]%
  }{%
    \end{proof}
  }

  \newtheorem*{TheoremA}{Theorem~A}
\begin{document}
%
%
%
  \newcommand{\Wisdom}[3]{%
    \begingroup
      #3
    \endgroup
  }
  \newcommand{\isom}{\cong}
  \newcommand{\updownarrows}{{\uparrow\downarrow}}
  \newcommand{\Tensor}[3]{#2\otimes_{#1}#3}
  \newcommand{\Homs}[3]{\Hom_{#1}(#2,#3)}
  \newcommand{\retraktarrows}{\mathrel{\mathpalette\retrakthelp{}}}
  \newcommand{\RetractDiagram}[2]{%
    \ifXYpic
    \xymatrix{
      {#1}
      \ar@{->>}@<0.75mm>[r]
    &
      {#2}
      \ar@{ >->}@<0.75mm>[l]
    }
    \else
      #1 \retractarrows #2      
    \fi
  }
  \makeatletter
  \newcommand{\retrakthelp}[1]{\vcenter{\m@th\hbox{\ooalign{\raise3pt
      \hbox{$#1\twoheadrightarrow$}\crcr$#1\leftarrowtail$}}}}
  \makeatother
  \newcommand{\mapcolon}{\colon\thinspace}
  \newcommand{\suchthatspace}{\thickspace}
  \newcommand{\suchthatvrule}{\suchthatspace\vrule\suchthatspace}
  \newcommand{\ssuchthatvrule}{\suchthatspace|\suchthatspace}
  \newcommand{\bsuchthatvrule}{\suchthatspace\big|\suchthatspace}
  \newcommand{\bbsuchthatvrule}{\suchthatspace\bigg|\suchthatspace}
%
%
  \newcommand{\FP}[1]{FP$_{#1}$}
  \newcommand{\F}[1]{F$_{#1}$}
  \newcommand{\CP}[1]{CP$_{#1}$}
  \newcommand{\type}{m}
  \newcommand{\Rank}{n}
  \newcommand{\vardim}{l}
  \newcommand{\ldec}[1]{\downarrow\vcenter{\rlap{$\scriptstyle #1$}}}
  \newcommand{\variable}{t}
  \newcommand{\CardOf}[1]{|#1|}
  \newcommand{\SetOf}[1]{\left\{#1\right\}}
  \newcommand{\sSetOf}[1]{\{#1\}}
  \newcommand{\bSetOf}[1]{\bigl\{#1\bigr\}}
  \newcommand{\bbSetOf}[1]{\biggl\{#1\biggr\}}
  \newcommand{\TupelOf}[1]{\left(#1\right)}
  \newcommand{\sTupelOf}[1]{(#1)}
  \newcommand{\bTupelOf}[1]{\bigl(#1\bigr)}
  \newcommand{\bbTupelOf}[1]{\biggl(#1\biggr)}
  \newcommand{\FamOf}[2]{\TupelOf{#1}_{{#2}}}
  \newcommand{\sFamOf}[2]{\sTupelOf{#1}_{{#2}}}
  \newcommand{\bFamOf}[2]{\bTupelOf{#1}_{{#2}}}
  \newcommand{\bbFamOf}[2]{\bbTupelOf{#1}_{{#2}}}
  \newcommand{\RankOf}[1]{\rank_{\ZZZ}#1}
  \newcommand{\rank}{\operatorname{rank}}
  \newcommand{\GroupGen}[1]{\left\langle #1 \right\rangle}
  \newcommand{\Closure}[1]{\overline{#1}}
  \newcommand{\Union}[2]{\bigcup_{#1}#2}
  \newcommand{\Intersection}[2]{\bigcap_{#1}#2}
  \newcommand{\CrossProd}[2]{\bigtimes_{#1}#2}
  \newcommand{\Sum}[2]{\sum_{#1}#2}
  \newcommand{\DisjointUnion}[2]{\bigudot_{#1}#2}
  \newcommand{\Ev}[2]{#1\left(#2\right)}
  \newcommand{\Kern}[1]{\ker #1}
  \newcommand{\GrHom}[2]{\Hom(#1,#2)}
  \newcommand{\Dual}[1]{#1^*}
  \newcommand{\PreIm}[1]{#1^{-1}}
  \newcommand{\PreImOf}[2]{\PreIm{#1}(#2)}
  \newcommand{\Restr}[2]{#1\mathord{\mid}_{#2}}
  \newcommand{\Pair}[2]{\left(#1,#2\right)}
  \newcommand{\compl}[1]{#1^{\mathrm{c}}}
  \newcommand{\bns}{\Sigma}
  \newcommand{\Bns}[1]{\bns^{#1}}
  \newcommand{\BNS}[2]{\Bns{#1}(#2)}
  \newcommand{\BNSc}[2]{\compl{\BNS{#1}{#2}}}
  \newcommand{\hght}{\Hight}
  \newcommand{\homo}{\chi}
  \newcommand{\homoat}[1]{\homo(#1)}
  \newcommand{\Height}[1]{\hght_{#1}}
  \newcommand{\HeightAt}[2]{\Height{#1}(#2)}
  \newcommand{\KomplexPoint}{y}
  \newcommand{\Komplex}{Y}
  \newcommand{\conv}{\operatorname{conv}}
  \newcommand{\Conv}[2]{\conv_{#1}(#2)}
  \newcommand{\ConvSub}[3]{\Conv{#1,#2}{#3}}
  \newcommand{\Koeffizient}{\mu}
  \newcommand{\KoeffizientAt}[1]{\Koeffizient_{#1}}
  \newcommand{\KoeffizientAtInd}[2]{\Koeffizient_{#2,#1}}
  \newcommand{\SigmaSet}{\Psi}
  \newcommand{\RootSubSet}{\Phi'}
  \newcommand{\Kegel}[1]{C_{#1}}
  \newcommand{\ConeSpannedBy}[1]{\overline{#1}}
  \newcommand{\UnitsOf}[1]{#1^*}
  \newcommand{\UnitSquaresOf}[1]{{\UnitsOf{#1}}^2}
  \newcommand{\EquivClass}[1]{\left[#1\right]}
  \newcommand{\Boundary}[1]{\partial #1}
  \newcommand{\Star}[1]{\mathrm{St}\left(#1\right)}
  \newcommand{\LeftPart}{L}
  \newcommand{\RightPart}{R}
  \newcommand{\Group}{G}
  \newcommand{\subgroup}{\leq}
  \newcommand{\SubGroup}{N}
  \newcommand{\FaktorGroup}{H}
  \newcommand{\GroupAt}[1]{\Group_{#1}}
  \newcommand{\AdeleGroup}{\boldsymbol{\Group}}
  \newcommand{\GroupScheme}{\LAG}
  \newcommand{\GroupElement}{g}
  \newcommand{\normalsubgroup}{\trianglelefteq}
  \newcommand{\Ring}{R}
  \newcommand{\KGOne}[1]{K(#1,1)}
  \newcommand{\Field}{K}
  \newcommand{\AnyField}{\Field'}
  \newcommand{\AnyRing}{R}
  \newcommand{\FieldAt}[1]{\Field_{#1}}
  \newcommand{\FieldEl}{x}
  \newcommand{\field}{k}
  \newcommand{\function}{f}
  \newcommand{\Set}{S}
  \newcommand{\OkaSet}{\Oka_{\Set}}
  \newcommand{\OkaAt}[1]{\Oka_{#1}}
  \newcommand{\maxidealat}[1]{\mathfrak{m}_{{#1}}}
  \newcommand{\locunifat}[1]{\pi_{{#1}}}
  \newcommand{\fieldat}[1]{\field_{{#1}}}
  \newcommand{\degat}[1]{d_{{#1}}}
  \newcommand{\Places}{P}
  \newcommand{\place}{p}
  \newcommand{\AdeleTo}[1]{\Adele_{{#1}}}
  \newcommand{\Adele}{\boldsymbol{A}}
  \newcommand{\adel}{\boldsymbol{a}}
  \newcommand{\RingEl}{t}
  \newcommand{\Over}[2]{#1(#2)}
  \newcommand{\LAG}{\mathcal{G}}
  \newcommand{\LAGOver}[1]{\Over{\LAG}{#1}}
  \newcommand{\Mult}{{\mathfrak{Mult}}}
  \newcommand{\Add}{{\mathfrak{Add}}}
  \newcommand{\Characters}[1]{\widehat{#1}}
  \newcommand{\CharactersOver}[2]{\Over{\Characters{#1}}{#2}}
  \newcommand{\character}{\chi}
  \newcommand{\CharacterInd}[1]{{\character}_{#1}}
  \newcommand{\elem}{g}
  \newcommand{\Aufl}{\mathcal{B}}
  \newcommand{\Chev}{\mathcal{G}}
  \newcommand{\TorDim}{n}
  \newcommand{\RkOf}[1]{\rk #1}
  \newcommand{\ChevTo}[1]{\Chev_{#1}}
  \newcommand{\ChevOver}[1]{\Over{\Chev}{#1}}
  \newcommand{\ChevOverTo}[2]{\Over{\ChevTo{#2}}{#1}}
  \newcommand{\ChevAt}[1]{\Chev_{#1}}
  \newcommand{\ChevOverAt}[2]{\Over{\ChevAt{#2}}{#1}}
  \newcommand{\Borel}{\mathcal{B}}
  \newcommand{\BorelOver}[1]{\Over{\Borel}{#1}}
  \newcommand{\Abels}{\mathcal{A}}
  \newcommand{\AbelsOver}[1]{\Over{\Abels}{#1}}
  \newcommand{\GenAbels}{A}
  \newcommand{\SubTorus}{\Torus'}
  \newcommand{\SubTorusOver}[1]{\Over{\SubTorus}{#1}}
  \newcommand{\GenSubTorus}{T'}
  \newcommand{\Unipot}{U}
  \newcommand{\Uni}{\mathcal{U}}
  \newcommand{\UniOver}[1]{\Over{\Uni}{#1}}
  \newcommand{\Torus}{\mathcal{T}}
  \newcommand{\TorusOver}[1]{\Over{\Torus}{#1}}
  \newcommand{\TorusElem}{\tau}
  \newcommand{\TorusElemAt}[1]{\TorusElem_{#1}}
  \newcommand{\RootSystem}{\Phi}
  \newcommand{\RootBasis}{\Delta^{-}}
  \newcommand{\PositiveRoots}{\RootSystem^{+}}
  \newcommand{\NegativeRoots}{\RootSystem^{-}}
  \newcommand{\NegativeRootsAt}[1]{\NegativeRoots_{#1}}
  \newcommand{\NegativeRootsUnion}{\boldsymbol{\NegativeRoots}}
  \newcommand{\NegativeRootsRestricted}{\widetilde{\NegativeRootsUnion}}
  \newcommand{\BaseRootsAt}[1]{\RootBasis_{#1}}
  \newcommand{\BaseRootsUnion}{\boldsymbol{\RootBasis}}
  \newcommand{\BaseRootsRestricted}{\widetilde{\BaseRootsUnion}}
  \newcommand{\NegativeRoot}{\root^{-}}
  \newcommand{\AffRoot}{\alpha}
  \newcommand{\AffRootAt}[1]{\AffRoot^{(#1)}}
  \newcommand{\AffRootInd}[1]{\AffRoot_{#1}}
  \newcommand{\AffRootAtInd}[2]{\AffRootAt{#1}_{#2}}
  \newcommand{\AffHight}{\Hight}
  \newcommand{\AffHightAt}[1]{\AffHight_{#1}}
  \newcommand{\Build}{X}
  \newcommand{\SphBuild}{\Delta}
  \newcommand{\SphBuildOver}[1]{\Delta\left(#1\right)}
  \newcommand{\AffBuild}{\Build}
  \newcommand{\AffBuildAt}[1]{\AffBuild_{#1}}
  \newcommand{\AffBuildOverAt}[2]{\Over{\AffBuildAt{#2}}{#1}}
  \newcommand{\InfBuild}{\SphBuild}
  \newcommand{\InfBuildAt}[1]{\InfBuild_{#1}}
  \newcommand{\InfBuildTo}[1]{{#1}^\infty}
  \newcommand{\InfRoot}{\AffRoot^{\infty}}
  \newcommand{\InfRootAt}[1]{\InfRoot_{#1}}
  \newcommand{\InfRootAtInd}[2]{{\InfRootAt{#1}}^{#2}}
  \newcommand{\AffApp}{\Sigma}
  \newcommand{\AffAppAt}[1]{\AffApp_{#1}}
  \newcommand{\AffAppAtInd}[2]{\AffAppAt{#1}^{#2}}
  \newcommand{\AffAppAtOf}[2]{\AffAppAt{#1}\left(#2\right)}
  \newcommand{\InfApp}{\Sigma^\infty}
  \newcommand{\InfAppAt}[1]{\InfApp_{#1}}
  \newcommand{\InfAppAtInd}[2]{{\AffAppAtInd{#1}{#2}}^\infty}
  \newcommand{\Chamber}{C}
  \newcommand{\ChamberInd}[1]{\Chamber_{#1}}
  \newcommand{\chamber}{c}
  \newcommand{\altchamber}{\tilde{c}}
  \newcommand{\AffChamber}{\chamber}
  \newcommand{\AffChamberInd}[1]{\AffChamber_{#1}}
  \newcommand{\AltChamber}{\altchamber}
  \newcommand{\AltChamberInd}[1]{\AffChamber_{#1}}
  \newcommand{\InfChamber}{\Chamber}
  \newcommand{\InfChamberAt}[1]{\InfChamber_{#1}}
  \newcommand{\InfChamberAtInd}[2]{\InfChamberAt{#1}^{#2}}
  \newcommand{\Sector}{S}
  \newcommand{\SectorAt}[1]{\Sector_{#1}}
  \newcommand{\SectorAtInd}[2]{\Sector_{#1}^{#2}}
  \newcommand{\gallery}{g}
  \newcommand{\Gallery}[1]{\left(#1\right)}
  \newcommand{\Panel}{\pi}
  \newcommand{\Autom}{\gamma}
  \newcommand{\AutomInd}[1]{\Autom_{#1}}
  \newcommand{\AutomIndEv}[2]{\Ev{\AutomInd{#1}}{#2}}
  \newcommand{\AutOf}[1]{\Aut\left(#1\right)}
  \newcommand{\IntersectionOne}{D}
  \newcommand{\Coordinates}{\varphi}
  \newcommand{\CoordinatesAt}[1]{\Coordinates_{#1}}
  \newcommand{\CoordinatesSum}{\boldsymbol{\Coordinates}}
  \newcommand{\Pro}{\pi}
  \newcommand{\ProAt}[1]{\Pro_{#1}}
  \newcommand{\RootGroup}{\U}
  \newcommand{\RootGroupTo}[1]{\RootGroup_{#1}}
  \newcommand{\RootGroupAtTo}[2]{\RootGroup_{#1,#2}}
  \newcommand{\AltAffRoot}{\beta}
  \newcommand{\AltAffRootInd}[1]{\AltAffRoot_{#1}}
  \newcommand{\RootGroupProdFromTo}[2]{\RootGroup^{(#1,#2)}}
  \newcommand{\RootGroupElement}{u}
  \newcommand{\RootGroupElementTo}[1]{\RootGroupElement_{#1}}
  \newcommand{\RootGroupElementUpTo}[1]{\RootGroupElement^{(#1)}}
  \newcommand{\RootIso}[1]{}
  \newcommand{\RootHom}[1]{}
  \newcommand{\RootHomEv}[2]{\RootHom{#1}\letf(#2\right)}
  \newcommand{\RootIsoEv}[2]{\RootIso{#1}\letf(#2\right)}
  \newcommand{\AffBuildProd}{\boldsymbol{\AffBuild}}
  \newcommand{\AffAppProd}{\boldsymbol{\AffApp}}
  \newcommand{\ProProd}{\boldsymbol{\Pro}}
  \newcommand{\Simplex}{\sigma}
  \newcommand{\SimplexAt}[1]{\Simplex_{#1}}
  \newcommand{\SimplexProd}{\boldsymbol{\Simplex}}
  \newcommand{\ProSum}{\boldsymbol{\xi}}
  \newcommand{\Kernel}{\boldsymbol{H}}
  \newcommand{\KernDual}{\Dual{\Kernel}}
  \newcommand{\SubKernel}{\Kernel_{\GenSubTorus}}
  \newcommand{\SubDual}{\Dual{\Kernel}_{\GenSubTorus}}
  \newcommand{\SubTransponiert}{\Transponiert{\SubKernel}}
  \newcommand{\Transponiert}[1]{#1^\tau}
  \newcommand{\BorelSpace}{\boldsymbol{Y}}
  \newcommand{\linform}{\SheetLin}
  \newcommand{\Tree}{T}
  \newcommand{\TreePoint}{\tau}
  \newcommand{\TreeInd}[1]{\Tree_{#1}}
  \newcommand{\TreeAt}[1]{\Tree_{#1}}
  \newcommand{\TreePointInd}[1]{\TreePoint_{#1}}
  \newcommand{\TreePointAt}[1]{\TreePoint_{#1}}
  \newcommand{\TreeProd}{\boldsymbol{\Tree}}
  \newcommand{\TreeProdTupel}{\boldsymbol{\TreePoint}}
  \newcommand{\TreeProdSub}[1]{\Sub{\TreeProd}{#1}}
  \newcommand{\sTreeProdSub}[1]{\sSub{\TreeProd}{#1}}
  \newcommand{\bTreeProdSub}[1]{\bSub{\TreeProd}{#1}}
  \newcommand{\bbTreeProdSub}[1]{\bbSub{\TreeProd}{#1}}
  \newcommand{\TreeHight}{\Hight^*}
  \newcommand{\TreeHightAt}[1]{\TreeHight_{#1}}
  \newcommand{\TreeHightSum}{\boldsymbol{\TreeHight}}
  \newcommand{\TreeSection}{\sigma}
  \newcommand{\TreeSectionEv}[1]{\TreeSection(#1)}
  \newcommand{\RestrTreeSection}{\sigma}
  \newcommand{\RestrTreeSectionEv}[1]{\TreeSection(#1)}
  \newcommand{\RestrTreeCorrection}{\rho}
  \newcommand{\RestrTreeCorrectionEv}[1]{\RestrTreeCorrection(#1)}
  \newcommand{\TreeProdPoint}{\boldsymbol{t}}
  \newcommand{\RestrTreeProd}{\boldsymbol{\Tree'}}
  \newcommand{\RestrTreeProdPoint}{\boldsymbol{\TreePoint'}}
  \newcommand{\TheSetOfLines}{\mathcal{L}}
  \newcommand{\AffAppPoint}{x}
  \newcommand{\AffAppPointAt}[1]{\AffAppPoint_{#1}}
  \newcommand{\AffAppPointAtInd}[2]{\AffAppPointAt{#1}^{#2}}
  \newcommand{\AffAppProdTupel}{\boldsymbol{\AffAppPoint}}
  \newcommand{\AffBuildPoint}{x}
  \newcommand{\AffBuildPointAt}[1]{\AffBuildPoint_{#1}}
  \newcommand{\AffBuildPointAtInd}[2]{\AffBuildPointAt{#1}^{#2}}
  \newcommand{\AffBuildProdTupel}{\boldsymbol{\AffAppPoint}}
  \newcommand{\SubComplProdAffApp}{Y}
  \newcommand{\SubComplProdAffBuild}{\tilde{Y}}
  \newcommand{\Complex}{X}
  \newcommand{\ComplexTo}[1]{\Complex_{#1}}
  \newcommand{\DirSet}{D}
  \newcommand{\direl}{\alpha}
  \newcommand{\direla}{\alpha}
  \newcommand{\direlb}{\beta}
  \newcommand{\RedHomTo}[2]{\widetilde{\H}_{#1}\left(#2\right)}
  \newcommand{\sRedHomTo}[2]{\widetilde{\H}_{#1}(#2)}
  \newcommand{\bRedHomTo}[2]{\widetilde{\H}_{#1}\bigl(#2\bigr)}
  \newcommand{\bbRedHomTo}[2]{\widetilde{\H}_{#1}\biggl(#2\biggr)}
  \newcommand{\HomTo}[2]{\H_{#1}(#2)}
  \newcommand{\FundamentalGroup}[2]{\pi_{#1}(#2)}
  \newcommand{\Sl}{\SL}
  \newcommand{\SlOver}[1]{\Sl(#1)}
  \newcommand{\SlInd}[1]{\SL_{#1}}
  \newcommand{\SlIndOver}[2]{\Over{\SlInd{#1}}{#2}}
  \newcommand{\SUTriInd}[1]{B^0_{#1}}
  \newcommand{\SUTriIndOver}[2]{\Over{\SUTriInd{#1}}{#2}}
  \newcommand{\SDiagInd}[1]{D^0_{#1}}
  \newcommand{\SDiagIndOver}[2]{\Over{\SDiagInd{#1}}{#2}}
  \newcommand{\UniTriInd}[1]{U_{#1}}
  \newcommand{\GlInd}[1]{\GL_{#1}}
  \newcommand{\GlIndOver}[2]{\Over{\GlInd{#1}}{#2}}
  \newcommand{\UTriInd}[1]{B_{#1}}
  \newcommand{\UTriIndOver}[2]{\Over{\UTriInd{#1}}{#2}}
  \newcommand{\DiagInd}[1]{D_{#1}}
  \newcommand{\DiagIndOver}[2]{\Over{\DiagInd{#1}}{#2}}
  \newcommand{\SSUTriInd}[1]{B^\pm_{#1}}
  \newcommand{\SSUTriIndOver}[2]{\Over{\SSUTriInd{#1}}{#2}}
  \newcommand{\Diag}[2]{\left(\begin{array}{cc}
                                 #1 & 0 \\
                                 0 & #2
                               \end{array}\right)}
  \newcommand{\PSUTriIndOver}[2]{\PPP\SUTriIndOver{#1}{#2}}
  \newcommand{\PUTriIndOver}[2]{\PPP\UTriIndOver{#1}{#2}}
  \newcommand{\Sub}[2]{#1\left[#2\right]}
  \newcommand{\sSub}[2]{#1[#2]}
  \newcommand{\bSub}[2]{#1\bigl[#2\bigr]}
  \newcommand{\bbSub}[2]{#1\biggl[#2\biggr]}
  \newcommand{\Hight}{h}
  \newcommand{\HightInd}[1]{\Hight_{#1}}
  \newcommand{\HightSum}{\boldsymbol{\Hight}}
  \newcommand{\Line}{L}
  \newcommand{\LineInd}[1]{\Line_{#1}}
  \newcommand{\LineIndDir}[2]{\LineInd{#1}^{#2}}
  \newcommand{\LineProd}{\boldsymbol{\Line}}
  \newcommand{\LinePoint}{{\tau'}}
  \newcommand{\LinePointInd}[1]{\LinePoint_{#1}}
  \newcommand{\LineProdTupel}{\boldsymbol{\LinePoint}}
  \newcommand{\Retraktion}{\rho}
  \newcommand{\RetraktionAt}[1]{\Retraktion_{#1}}
  \newcommand{\RetraktionProd}{\boldsymbol{\Retraktion}}
  \newcommand{\Intervall}{I}
  \newcommand{\IntervallInd}[1]{\Intervall_{#1}}
  \newcommand{\IntervallProd}{\boldsymbol{\Intervall}}
  \newcommand{\SupIntervall}{J}
  \newcommand{\SupIntervallInd}[1]{\SupIntervall_{#1}}
  \newcommand{\SupIntervallProd}{\boldsymbol{\SupIntervall}}
  \newcommand{\Kompaktum}{\mathsf{C}}
  \newcommand{\KompaktDomain}{\Kompaktum'}
  \newcommand{\AltKompaktDomain}{\Kompaktum''}
  \newcommand{\SigmaComplOverTo}[2]{\Sigma^\mathrm{c}_{#1}(#2,\ZZZ)}
  \newcommand{\PrimeNumber}{p}
  \newcommand{\SheetCompl}{X}
  \newcommand{\Eukl}{\EEE}
  \newcommand{\EuKern}{\Eukl'}
  \newcommand{\SheetProj}{\pi}
  \newcommand{\SheetSet}{\SigmaSet}
  \newcommand{\SheetLin}{\lambda}
  \newcommand{\SheetLinEv}[1]{\SheetLin(#1)}
  \newcommand{\SheetConst}{c}
  \newcommand{\SheetConstAt}[1]{\SheetConst_{#1}}
  \newcommand{\SheetSubSet}{\SigmaSet'}
  \newcommand{\AffHalfSpace}{\AAA}
  \newcommand{\SheetPreim}{\SheetCompl'}
  \newcommand{\Sheet}{B}
  \newcommand{\SheetInd}[1]{\Sheet_{#1}}
  \newcommand{\SheetFiltrInd}[1]{W^{#1}}
  \newcommand{\Translationspace}{\TTT}
  \newcommand{\EuklPoint}{e}
  \newcommand{\SheetNewPiece}{\ConvexStep'}
  \newcommand{\ConvexStep}{V}
  \newcommand{\BuildPunkt}{x}
  \newcommand{\BuildPunktAt}[1]{\BuildPunkt^{(#1)}}
  \newcommand{\AppPunkt}{s}
  \newcommand{\StApp}{\Sigma_0}
  \newcommand{\StAppProd}{\overline{\StApp}}
  \newcommand{\StAppAt}[1]{\StApp^{(#1)}}
  \newcommand{\Transl}{\tau}
  \newcommand{\TranslAt}[1]{\Transl^{(#1)}}
  \newcommand{\retr}{\rho}
  \newcommand{\OkaTrans}{H}
  \newcommand{\Preimage}{\BuildProd_0}
  \newcommand{\TranslSpace}{\TTT}
  \newcommand{\Faktor}{\RRR^m}
  \newcommand{\SB}{\mathrm{B}^0}
  \newcommand{\Sph}[1]{\SSS(#1)}
  \newcommand{\nat}{n}
  \newcommand{\Point}{x}
  \newcommand{\cell}{\sigma}
  \newcommand{\celldim}{i}
  \newcommand{\Dir}{D}
  \newcommand{\dirind}{\alpha}
  \newcommand{\dirindzwei}{\beta}
  \newcommand{\AbGruppe}{A}
  \newcommand{\height}{h}
  \newcommand{\firstind}{1}
  \newcommand{\lastind}{\BaseDim}
  \newcommand{\LinRoot}{\aroot}
  \newcommand{\LinRootAt}[1]{\LinRoot^{(#1)}}
  \newcommand{\LinRootTo}[1]{\LinRoot_{#1}}
  \newcommand{\RootDist}{\delta}
  \newcommand{\Face}{\sigma}
  \newcommand{\FaceTo}[1]{\Face_{#1}}
  \newcommand{\ObDrei}{\UTri}
  \newcommand{\Base}{\EEE}
  \newcommand{\BaseDim}{n}
  \newcommand{\Hyp}{\mathcal{H}}
  \newcommand{\Ebene}{H}
  \newcommand{\Gew}{\XXX}
  \newcommand{\Sch}{\Sigma}
  \newcommand{\SchTo}[1]{\Sch_{#1}}
  \newcommand{\SchEins}{\SchTo{1}}
  \newcommand{\SchZwei}{\SchTo{2}}
  \newcommand{\SchPunkt}{\AppPunkt}
  \newcommand{\SchPunktTo}[1]{\SchPunkt_{#1}}
  \newcommand{\SchPunktEins}{\SchPunktTo{1}}
  \newcommand{\SchPunktZwei}{\SchPunktTo{2}}
  \newcommand{\Gerade}{g}
  \newcommand{\GeradeIn}[1]{\Gerade_{#1}}
  \newcommand{\TREE}{\mathcal{Y}}
  \newcommand{\DirEl}{\tau}
  \newcommand{\PosWurzelZu}[2]{\Wurzel^+_{#1}(#2)}
  \newcommand{\NegWurzelZu}[2]{\Wurzel^-_{#1}(#2)}
  \newcommand{\Directions}{\Delta}
  \newcommand{\BaseDirections}{\Directions^0}
  \newcommand{\BaseDirEl}{{\DirEl^0}}
  \newcommand{\BuildToTree}{\GewToTree}
  \newcommand{\BuildToTreeAt}[1]{\BuildToTree^{(#1)}}
  \newcommand{\TreeToBuild}{\TreeToGew}
  \newcommand{\TreeToBuildAt}[1]{\TreeToBuild^{(#1)}}
  \newcommand{\GewToTree}{\rho}
  \newcommand{\TreeToGew}{\varsigma}
  \newcommand{\TreeZu}[1]{\Tree_{#1}}
  \newcommand{\GewPunkt}{x}
  \newcommand{\GewPunktEins}{\GewPunkt_1}
  \newcommand{\GewPunktZwei}{\GewPunkt_2}
  \newcommand{\BasePunkt}{e}
  \newcommand{\BasePunktEins}{\BasePunkt_1}
  \newcommand{\BasePunktZwei}{\BasePunkt_2}
  \newcommand{\Span}{\operatorname{span}}
  \newcommand{\BaseSpan}[1]{\Span_{\Base}(#1)}
  \newcommand{\GewSpan}[1]{\Span_{\Gew}(#1)}
  \newcommand{\aequivalent}[3]{#2 \sim_{#1} #3}
  \newcommand{\LineZu}[1]{\RRR_{#1}}
  \newcommand{\LinePunkt}{l}
  \newcommand{\LinePunktEins}{\LinePunkt_1}
  \newcommand{\LinePunktZwei}{\LinePunkt_2}
  \newcommand{\BaseLine}{g}
  \newcommand{\coordinate}{x}
  \newcommand{\IntervallTo}[1]{\Intervall_{#1}}
  \newcommand{\IntervallFam}{\underline{\Intervall}}
  \newcommand{\GewIntervall}{\GewProd_{\IntervallFam}}
  \newcommand{\subcompl}{\subseteq}
  \newcommand{\IndexSch}{i}
  \newcommand{\MaxIndex}{m}
  \newcommand{\SupIndex}{j}
  \newcommand{\SubIndex}{i}
  \newcommand{\BotIndex}{\One}
  \newcommand{\One}{1}
  \newcommand{\TopIndex}{r}
  \newcommand{\LowReal}{a}
  \newcommand{\HighReal}{b}
  \newcommand{\TheLag}{l}
  \newcommand{\Index}{i}
  \newcommand{\IndexAt}[1]{\Index_{#1}}
  \newcommand{\LowInd}{r}
  \newcommand{\HighInd}{s}
  \newcommand{\IndexSet}{J}
  \newcommand{\Two}{2}
  \newcommand{\BlockTop}{q}
  \newcommand{\GroupTop}{k}
  \newcommand{\BlockInd}{i}
  \newcommand{\BlockSubInd}{i}
  \newcommand{\BlockSupInd}{j}
  \newcommand{\GroupSubInd}{\mu}
  \newcommand{\GroupSupInd}{\nu}
  \newcommand{\GroupInd}{\mu}
  \newcommand{\real}{t}
  \newcommand{\lowreal}{s}
  \newcommand{\RealAt}[1]{\real_{#1}}
  \newcommand{\RealAtInd}[2]{\RealAt{#1,#2}}
  \newcommand{\RealInd}[1]{\real_{#1}}
  \newcommand{\Exp}[1]{e^{{#1}}}
  \newcommand{\RootGroupAt}[1]{\RootGroup_{#1}}
  \newcommand{\AllBetween}[2]{[#1,#2]}
  \newcommand{\InBetween}[2]{(#1,#2)}
  \newcommand{\RootGroupEl}{u}
  \newcommand{\mouf}[1]{m(#1)}
  \newcommand{\MoufEl}{n}
  \newcommand{\panel}{\pi}
  \newcommand{\Apartment}{\Sigma}
  \newcommand{\ApartmentInd}[1]{\Apartment_{#1}}
  \renewcommand{\root}{\alpha}
  \newcommand{\roota}{\alpha}
  \newcommand{\rootb}{\beta}
  \newcommand{\rootc}{\gamma}
  \newcommand{\CW}{{\small CW}}
%
%
%
\newcommand{\CaseTwoAlpha}{%
  \setlength{\unitlength}{2.5mm}
  \begin{picture}(8,12)(0,0)
    \put(2,10){$\,\EuKern$}
    \put(5,10){$\AffHalfSpace$}
    \put(5,4){$\SheetNewPiece$}
    \thinlines\drawline(2,0)(2,12)
    \thicklines\drawline(8,0)(2,3)
    \thinlines\drawline(2,3)(0,4)(0,6)(2,7)
    \thicklines\drawline(2,7)(8,10)
  \end{picture}%
}%
\newcommand{\CaseOne}{%
  \setlength{\unitlength}{2.5mm}
  \begin{picture}(8,12)(0,0)
    \put(2,10){$\,\EuKern$}
    \put(5,10){$\AffHalfSpace$}
    \put(5,4){$\SheetNewPiece$}
    \thinlines\drawline(2,0)(2,12)
    \thinlines\drawline(0,0)(2,1)
    \thicklines\drawline(2,1)(8,4)(8,6)(2,9)
    \thinlines\drawline(2,9)(0,10)
  \end{picture}%
}%
\newcommand{\CaseTwoBeta}{%
  \setlength{\unitlength}{2.5mm}
  \begin{picture}(8,12)(0,0)
    \put(2,10){$\,\EuKern$}
    \put(5,10){$\AffHalfSpace$}
    \put(5,4){$\SheetNewPiece$}
    \thinlines\drawline(2,0)(2,12)
    \thinlines\drawline(0,2)(2,4)
    \thicklines\drawline(2,4)(6,8)(8,8)(8,0)
  \end{picture}%
}%
%
%
%
\title{Finiteness properties of soluble arithmetic groups\\over global function fields}%
\shorttitle{Finiteness properties of soluble arithmetic groups}

\author{Kai-Uwe Bux}%
\address{%
  Cornell University, Department of Mathemtics\\
  Malott Hall 310, Ithaca, NY 14853-4201, USA}
\gtemail{\href{mailto:bux_math_2004@kubux.net}{bux\_math\_2004@kubux.net}}
\asciiemail{bux_math_2004@kubux.net}
\urladdr{http://www.kubux.net}

\asciiabstract{Let G be a Chevalley group scheme and B<=G a Borel
  subgroup scheme, both defined over Z.  Let K be a global function
  field, S be a finite non-empty set of places over K, and O_S be the
  corresponding S-arithmetic ring.  Then, the S-arithmetic group
  B(O_S) is of type F_{|S|-1} but not of type FP_{|S|}.  Moreover one
  can derive lower and upper bounds for the geometric invariants
  \Sigma^m(B(O_S)). These are sharp if G has rank 1. For higher ranks,
  the estimates imply that normal subgroups of B(O_S) with abelian
  quotients, generically, satisfy strong finiteness conditions.}

\begin{abstract}
  Let $\Chev$ be a Chevalley group scheme
  and $\Borel\leq\Chev$ a Borel subgroup
  scheme, both defined over $\ZZZ$. Let $\Field$ be a global function field,
  $\Set$ be a finite non-empty set of places over $\Field$, and
  $\OkaSet$ be the corresponding $\Set$--arithmetic ring.
  Then, the $\Set$--arithmetic group
  $\BorelOver{\OkaSet}$ is of type \F{\CardOf{\Set}-1}
  but not of type \FP{\CardOf{\Set}}. Moreover one can derive lower and
  upper bounds for the geometric
  invariants $\BNS{\type}{\BorelOver{\OkaSet}}$. These are sharp
  if $\Chev$ has rank $1$. For higher ranks, the estimates imply
  that normal subgroups of $\BorelOver{\OkaSet}$ with abelian
  quotients, generically, satisfy strong finiteness conditions.
\end{abstract}
%
%
%

\maketitlepage

\noindent
Once upon a time, all finitely generated groups were finitely presented.
There were discontinuous subgroups of Lie groups or groups acting nicely
on beautiful geometries---one should think, for example, of finitely generated
Fuchsian groups. Then B\,H  Neumann gave the first example of
a finitely generated but infinitely related group in~\cite{Neumann:1937}
where he even showed that there are uncountably many $2$~generator
groups. Among these, of course, only countably many are finitely presented.
Hence finite presentability is a much stronger property than finite
generation.

More than twenty years later, working on decision problems,
G Baumslag, W\,W Boone, and B\,H Neumann
showed that even finitely generated subgroups of finitely presented
groups need not be finitely presented \cite{Baumslag.Boone.Neumann:1959}.

However, it took almost another twenty years until finitely generated infinitely related
groups were ``observed in nature''. In \cite{Stuhler:1976}, U Stuhler
proved that the groups $\SlIndOver{2}{\field[\variable,\variable^{-1}]}$,
where $\field$ is a finite field, are finitely generated but not finitely presented.
He extended these results in
\cite{Stuhler:1980} constructing series of groups with increasing
finiteness properties.
In \cite{Strebel:1984}, R Strebel gives a historical and systematic survey
with focus on soluble groups.

Because of the topological background of finite generation and
finite presentability, there are two generalizations to higher
dimensions: one based on homotopy groups, the other based on homology.
A group $\Group$ is \notion{of type \F{\type}} if
there is an Eilenberg--MacLane space $\KGOne{\Group}$ with
finite $\type$--skeleton. $\Group$ is \notion{of type \FP{\type}}
if the trivial $\ZZZ\Group$--module $\ZZZ$ admits a projective
resolution that is finitely generated in dimensions $\leq\type$.
This is a homological variant of the homotopical finiteness
property \F{\type}. The definition of type \F{\type} was given
by  C\,T\,C Wall in \cite{Wall:1965}. It is convenient to define
the finiteness length of a group to be the largest dimension
$\type$ for which a group is of type \F{\type}. 

A group is of type \F{1} if and only if it is finitely generated.
Moreover, type \F{1} and type \FP{1} are equivalent notions.
A group is of type \F{2} if and only if it is finitely presented.
M Bestvina and N Brady \cite{Bestvina.Brady:1997} have given
an example of a group of type \FP{2} that is not
finitely presented. However, this is the only way in which homotopical
and homological finiteness properties differ:
\begin{quote}
  For $\type\geq2$, a group $\Group$ is of type \F{\type} if and only if it is
  finitely presented and of type \FP{\type}.
\end{quote}

Finiteness properties are still somewhat mysterious. Theorems relating finiteness
properties in a transparent way to other, more group theoretic properties are
in short supply. For special classes of groups, however, the situation is
better. Eg, for metabelian groups, the
Bieri--Strebel theory of geometric invariants leads to nice conjectures
which are confirmed by a lot of examples and partial results.
Finite presentability is well understood within this context. Eg, these
groups are finitely presented if and only if they are of type \FP{2} \cite{Bieri.Strebel:1980}.

In a situation like this, the best one can hope for is to get a feeling for
the relationship between finiteness conditions and group structure
within certain classes of groups. We will consider a class of
$\Set$--arithmetic groups. These matrix groups are given by means of
an algebraic group scheme $\GroupScheme$ and a set $\Set$ of primes
over a global field $\Field$ which determines an $\Set$--arithmetic
ring $\OkaSet\subseteq\Field$. These two parameters can be
varied independently, and one would like to know how finiteness properties vary
with them. Moreover, these groups are natural generalizations
of lattices in Lie groups, for which finiteness properties often have a
more direct geometric interpretation. For all these reasons, a lot of research has
already been done on finiteness properties of $\Set$--arithmetic groups.

The theory of $\Set$--arithmetic groups is dominated by two fundamental
distinctions. The field $\Field$ can be a global number field or a global
function field. With respect to the group scheme, there are two extremes
the first of which is given by reductive groups,
eg, $\GlInd{\Rank}$ or $\SlInd{\Rank}$, which is even a Chevalley group.
Soluble groups, eg, groups of upper triangular
matrices form the other extreme.
Let us recall the most important results:
\begin{description}
  \item[$\GroupScheme$ reductive:]
    \begin{description}
      \item[]
      \item[$\Field$ number field:]
        $\Over{\GroupScheme}{\OkaSet}$ is of \notion{type \F{\infty}},
        ie, of type \F{\type} for all
        $\type \in \NNN$ \cite{Borel.Serre:1976}.
      \item[$\Field$ function field:]
        Finite generation and finite presentability are completely
        understood \cite{Behr:1992}.

        There are series of examples that support the conjecture
        that the finiteness length 
        grows with
        $\CardOf{\Set}$ and the rank of $\Chev$. The most important results are:
        \begin{itemize}
          \item
            $\Over{\SlInd{2}}{\OkaSet}$ is of type \F{\CardOf{\Set}-1} but
            not of type \FP{\CardOf{\Set}} \cite{Stuhler:1980}.
          \item
            If $\GroupScheme$ is a Chevalley group of rank $\TorDim$
            not of exceptional type, then
            $\Over{\GroupScheme}{\field[t]}$ is of type \F{\TorDim-1} but not
            of type \FP{\TorDim} provided the finite field~$\field$ is
            large enough \cite[Corollary~20, page~113]{Abramenko:1996}.
        \end{itemize}
    \end{description}
  \item[$\Aufl$ soluble:]
    \begin{description}
      \item[]
      \item[$\Field$ number field:]
        Finite generation and finite presentability are completely
        understood. Finite presentability is treated in~\cite{Abels:1987}.

        In~\cite[Theorem~3.1]{Tiemeier:1997}, a Hasse principle is
        derived:
        $\Over{\Aufl}{\OkaSet}$ is of type \FP{\type} if and only if
        for each place $\place\in\Set$, the group
        $\Over{\Aufl}{\OkaAt{\place}}$ satisfies the
        compactness property~\CP{\type}, which is defined
        in~\cite{Abels.Tiemeier:1997}.
        If $\Aufl$ is a Borel subgroup of a Chevalley group,
        $\Over{\Aufl}{\OkaSet}$ is of type
        \F{\infty} \cite[Corollary~4.5]{Tiemeier:1997}.

        Beyond these, there are some series of examples, eg,
        in~\cite{Abels.Brown:1987}.
      \item[$\Field$ function field:]
        No Hasse principle holds in this case. This follows already from
        the series of metabelian groups that is examined using
        Bieri--Strebel theory in~\cite{Bux:1997}.
        In this article, we generalize those results to group
        schemes of higher rank.
    \end{description}
\end{description}
Our main result is the following:

\begin{TheoremA}
  Let $\Chev$ be a Chevalley group, $\Borel\leq\Chev$ a Borel subgroup,
  $\Field$ a global function field, $\Set$ a non-empty set of places
  over $\Field$, and $\OkaSet$ the corresponding $\Set$--arithmetic ring.
  Then $\BorelOver{\OkaSet}$ is of type \F{\CardOf{\Set}-1} but not
  of type \FP{\CardOf{\Set}}.
\end{TheoremA}

We will define notions and fix notations in the first two sections.
Then we will deal with the rank--$1$-case. Sections~\ref{BeginProof}
to~\ref{EndProof} contain the proof of Theorem~A: in
Theorem~\ref{ObereAbschaetzung} the upper bound is established
whereas the lower bound is given in Theorem~\ref{UntereAbschaetzung}.
The last section presents Theorem~\ref{Sigma}, which provides bounds
for the ``geometric invariants''. Definitions and a bit of motivation
will be given at the beginning of Section~\ref{Sec:Sigma}.

%
%
This paper grew out of my PhD thesis,
which I wrote under the supervision
of Prof Robert~Bieri. I would like to thank
him for his support and
encouragement. I also would like to thank the
referee for very carefully reading the paper and
suggesting numerous improvements.

\section{Preliminaries on ad\`{e}les and unipotent groups}
  General references about global fields and ad\`{e}les
  are~\cite{Cassels.Froehlich:1967} or~\cite{Weil:1973}.
  In this paper
  \begin{notationlist}
    \item
      $\Field$ is a global function field. Its elements are called
      \notion{functions}. Let
    \item
      $\Places$ denote the set of all places of $\Field$. We regard a place
      as a normalized discrete valuation
      $\place \mapcolon \Field \rightarrow \ZZZ \cup \SetOf{ \infty }$.
      For each place $\place$, there is a \notion{local function field}
    \item
      $\FieldAt{\place}$, which is the completion of $\Field$ at $\place$.
      This is a topological field.
      Extending $\place$ continuously, we obtain a normalized discrete
      valuation on $\FieldAt{\place}$, which we also denote by $\place$.
      The subring of functions \notion{holomorphic at $\place$} is denoted by
    \item
      $\OkaAt{\place} := 
      \SetOf{ \function \in \FieldAt{\place} \suchthatvrule \place(\function) \geq 0 }$.
      This ring is a compact open subspace of $\FieldAt{\place}$. Moreover, it
      is a local ring with maximal ideal
    \item
      $\maxidealat{\place} = 
        \SetOf{ \function \in \FieldAt{\place} \suchthatvrule \place(\function) > 0 }$.
      The \notion{constant} functions, ie, the functions in $\Field$ that are
      holomorphic everywhere form a finite subfield
    \item
      $\field$.
      The \notion{residue field}
    \item
      $\fieldat{\place} := \OkaAt{\place} \rmod \maxidealat{\place}$
      is a finite extension of the field $\field$ of degree
    \item
      $\degat{\place} := [ \fieldat{\place} \mapcolon \field ]$.
      We define a norm on $\FieldAt{\place}$ by
    \item
      $|\function|_\place := \Exp{-\degat{\place} \place(\function)}$,
      which is proportional to the \xnotion{modulus} of
      $\FieldAt{\place}$.
  \end{notationlist}
  Ad\`{e}les provide a formalism to dealt with all places simultaneously.
  For a finite non-empty set $\Set\subseteq\Places$ of places,
  \begin{notationlist}
    \item
      $\OkaSet := \SetOf{ \function \in \Field 
        \suchthatvrule \place(\function) \geq 0 \;\;\forall 
        \place \in \Places \setminus \Set }$
      is the ring of functions holomorphic outside $\Set$ and
    \item
      $\AdeleTo{\Set} :=
        \bigtimes_{\place \in \Set} \FieldAt{\place}
        \times
        \bigtimes_{\place \not\in \Set} \OkaAt{\place}$
      is the topological ring of $\Set$--ad\`{e}les.
  \end{notationlist}
  The \notion{ring of ad\`{e}les} of $\Field$ is the direct limit
  \begin{notationlist}
    \item
      $\Adele := \directlim_{\Set} \AdeleTo{\Set}$ of topological rings. For an
      \notion{ad\`{e}le} $\adel = \FamOf{\function_\place}{\place \in \Places}
       \in \Adele$, let
    \item
      $|\adel| := \prod_\place |\function_\place|_\place$
      be the \notion{id\`{e}lic norm}.
  \end{notationlist}
  There is an inclusion $\Field\subseteq\Adele$ because every function
  in $\Field$ is holomorphic at almost all places. This way,
  $\Field$ is a discrete subring of $\Adele$, and we have
  \[
    \OkaSet = \AdeleTo{\Set} \cap \Field
    .
  \]
  Since we are working with topological rings,
  we regard a linear algebraic group $\LAG$ defined over a
  commutative ring $\Ring$
  as a functor from the category of
  (topological) commutative $\Ring$--algebras into the category
  of (topological) groups.
  
  A fundamental theorem
  says that the quotient $\Adele \rmod \Field$ is
  compact \cite[Theorem~2, page~64]{Weil:1973}.
  We need a slightly more general statement.
  \begin{Lemma}\label{eins:kokompakt}
    Let $\Unipot$ be a unipotent linear algebraic group defined over
    $\Field$. Then
    $\Over{\Unipot}{\Field}$ is a discrete subgroup of
    $\Over{\Unipot}{\Adele}$ and the quotients
    $\Over{\Unipot}{\Adele} \rmod \Over{\Unipot}{\Field}$ and
    $\Over{\Unipot}{\Field} \lmod \Over{\Unipot}{\Adele}$ are compact.
  \end{Lemma}
  \begin{proof}
    Let $\Characters{\Unipot}$ denote the \notion{group of characters} of
    $\Unipot$, ie, the linear algebraic group of $\Field$--morphisms
    from $\Unipot$ to $\Mult:=\GlInd{1}$. Put
    \[
      \Over{\Unipot}{\Adele}^\circ :=
      \bSetOf{ 
              \elem \in \Over{\Unipot}{\Adele}
              \bsuchthatvrule
              |\character( \elem )| = 1
              \text{\ for all\ } \character \in \CharactersOver{\Unipot}{\Field} 
            }
      .
    \]
    According to ~\cite[Theorem~4.8]{Borel:1991}, $\Unipot$ can be
    triagonalized over $\Field$. Hence, \cite[Satz~3]{Behr:1969}
    implies that the quotient
    $\Over{\Unipot}{\Adele}^\circ \rmod \Over{\Unipot}{\Field}$
    is compact. Therefore, it suffices to prove that
    $\Over{\Unipot}{\Adele}^\circ = \Over{\Unipot}{\Adele}$. This,
    however, follows from the fact that unipotent groups do not admit
    non-trivial characters, since homomorphisms of affine algebraic
    groups preserve unipotency \cite[Theorem~4.4]{Borel:1991}.
  \end{proof}

\section{Chevalley groups and the associated buildings}
  We fix a \notion{Chevalley group}
  \begin{notationlist}
    \item
      $\Chev$, ie, a semisimple linear algebraic group defined over
      $\ZZZ$. How to build such a group scheme is described, for example,
      in~\cite{Chevalley:1960}, \cite{Steinberg:1967}, and~\cite{Abramenko:1994}.
      A standard references for the following facts are \cite{Bruhat.Tits:1972}
      and \cite{Bruhat.Tits:1984}.
      
      Every Chevalley group comes with a \xnotion{root system}
    \item
      $\RootSystem$. For any commutative unitary ring, we denote the
      group of $\AnyRing$--points of $\Chev$ by
    \item
      $\ChevOver{\AnyRing}$. If $\AnyField$ is a field, the group
      $\ChevOver{\AnyField}$
      acts \xnotion{strong-transitively}
      on the associated \xnotion{spherical building}
    \item
      $\SphBuildOver{\AnyField}$ of type $\RootSystem$. Given a place $\place$ on $\AnyField$,
      there is an affine building
    \item
      $\AffBuildOverAt{\AnyField}{\place}$ associated to
      and acted upon by $\ChevOver{\AnyField}$. Affine buildings have
      apartments that are Euclidean Coxeter complexes. Using the Euclidean
      metric on apartments, the induced path-metric on an affine building
      is {\small CAT$(0)$}. Moreover, apartments in affine buildings
      can be characterized in terms of this metric: a subspace is an
      apartment (in the complete system of apartments)
      if and only if it is a
      maximal flat subspace, ie, an isometrically embedded Euclidean
      space of maximal dimension.
      An excellent source and reference for the theory of affine buildings
      is \cite[Chapter~VI]{Brown:1989}.

      The building
      $\SphBuildOver{\AnyField}$ can be viewed in a natural way as the
      \xnotion{building at infinity} of $\AffBuildOverAt{\AnyField}{\place}$.
      If $\AnyField$ is complete with respect to $\place$, then
      the system of \xnotion{apartments} in $\SphBuildOver{\AnyField}$ induces the complete
      system of apartments in $\AffBuildOverAt{\AnyField}{\place}$.
      Then we have a
      $1$--$1$--correspondence of spherical and affine apartments.
      In this case, we say that an affine apartment contains a chamber,
      a point, or a halfapartment at infinity if the corresponding
      spherical apartment does. Put
    \item
      $\SphBuild:=\SphBuildOver{\Field}$ for the fixed global function field $\Field$. Let
    \item
      $\AffBuildAt{\place} := \AffBuildOverAt{\FieldAt{\place}}{\place}$
      denote the affine building associated
      to~$\ChevOver{\FieldAt{\place}}$ whereas
    \item
      $\InfBuildAt{\place}$ will denote the
      corresponding spherical building at infinity.

      We fix a chain
      \[
        \Torus \subgroup \Borel \subgroup \Chev
      \]
      of group schemes defined over $\ZZZ$ such that $\TorusOver{\Field}$ is
      a maximal $\Field$--split torus in $\ChevOver{\Field}$ and $\BorelOver{\Field}$
      is a \notion{Borel subgroup}, ie, a maximal solvable $\Field$--subgroup
      in $\ChevOver{\Field}$. Then, there is a unique apartment
    \item
      $\AffAppAt{\place}$ in $\AffBuildAt{\place}$ that is stabilized by
      $\TorusOver{\FieldAt{\place}}$. We regard $\AffAppAt{\place}$ as
      the \notion{standard apartment}. The group $\TorusOver{\FieldAt{\place}}$ acts on
      $\AffAppAt{\place}$ as a maximum rank lattice of translations. Thus, the
      action of $\TorusOver{\FieldAt{\place}}$ on $\AffAppAt{\place}$ is
      \notion{cocompact}, ie, the quotient space of orbits of this action is
      compact.
      
      Let
    \item
      $\TorDim$ be the dimension of $\Torus$.
      This is by definition the \notion{rank} of $\Chev$. The building
      $\AffBuildAt{\place}$ is a piecewise Euclidean complex of
      dimension $\TorDim$, and $\ChevOver{\FieldAt{\place}}$ acts upon $\AffBuildAt{\place}$ by
      \notion{cell-permuting} isometries: every element of a cell stabilizer fixes
      the cell pointwise.
      Cell stabilizers are open and compact.
      Since $\FieldAt{\place}$ is complete with respect to $\place$, the group
      $\ChevOver{\FieldAt{\place}}$ acts \notion{strongly-transitively}
      on $\AffBuildOverAt{\FieldAt{\place}}{\place}$, ie, the action is
      transitive on the set of pairs $(\AffApp,\AffChamber)$ where
      $\AffApp$ is an apartment in
      $\AffBuildOverAt{\FieldAt{\place}}{\place}$ and
      $\AffChamber$ is a chamber in $\AffApp$.
 
      The group $\BorelOver{\FieldAt{\place}}$
      is the stabilizer of a chamber at infinity
    \item
      $\InfChamberAt{\place}$ in the standard apartment $\AffAppAt{\place}$. We call this
      the
      \notion{fundamental chamber} at infinity. It is
      represented by a parallelity class of \xnotion{sectors} in
      $\AffBuildAt{\place}$.
      There is a canonical projection $\Borel \rightarrow \Torus$, which
      turns the torus $\Torus$ into a retract of $\Borel$. Let
    \item
      $\Uni$ denote its kernel, which is called
      the \notion{unipotent part} of $\Borel$.
      The group $\UniOver{\FieldAt{\place}}$ not just stabilizes the
      fundamental chamber $\InfChamberAt{\place}$, it \notion{fixes}
      this chamber at infinity, ie, for each
      element in $\UniOver{\FieldAt{\place}}$, there is a sector
      representing $\InfChamberAt{\place}$ which is fixed pointwise by
      the chosen element. This follows from the way the affine building
      $\AffBuildAt{\place}$ and the action of $\ChevOver{\FieldAt{\place}}$
      on $\AffBuildAt{\place}$ are constructed:
      the group $\UniOver{\FieldAt{\place}}$
      turns out to be generated by \notion{root groups}
      all of whose elements actually fix a Euclidean half apartment
      in $\AffAppAt{\place}$ containing $\InfChamberAt{\place}$.
      The construction of the building is described in
      \cite[Section~7.4]{Bruhat.Tits:1972} and the property of
      root groups used here is spelled out in
      Proposition~(7.4.5) of said section.
      (These root groups are spherical and not to be confused with the
      affine root groups discussed in Section~\ref{Moufang}.)

      Since any element of $\UniOver{\FieldAt{\place}}$ fixes a
      sector representing $\InfChamberAt{\place}$, it cannot move
      chambers within the standard apartment $\AffAppAt{\place}$ at all:
      just consider a gallery in $\AffAppAt{\place}$
      from a moved chamber to a chamber in the fixed sector.
      Since, on the other hand,
      $\ChevOver{\FieldAt{\place}}$ acts strongly-transitively,
      $\UniOver{\FieldAt{\place}}$ acts on
      $\AffBuildAt{\place}$ with $\AffAppAt{\place}$ as a
      \xnotion{fundamental domain}.
      We thus obtain a projection map
    \item
      $\ProAt{\place} \mapcolon \AffBuildAt{\place} \rightarrow \AffAppAt{\place}$.

      In order to determine the finiteness properties of $\BorelOver{\OkaSet}$,
      we will study its action on the product
    \item
      $\AffBuildProd := \bigtimes_{\place \in \Set} \AffBuildAt{\place}$
      of affine buildings. The projections $\ProAt{\place}$ induce a map
    \item
      $\ProProd \mapcolon \AffBuildProd \rightarrow \AffAppProd :=
      \bigtimes_{\place \in \Set} \AffAppAt{\place}$ onto the product
      of standard apartments.
  \end{notationlist}

  \begin{Lemma}\label{Projektion:eigentlich}
    The map $\ProProd$ induces a proper map
    $\UniOver{\OkaSet} \lmod \AffBuildProd \rightarrow \AffAppProd$.
  \end{Lemma}
  \begin{proof}
    Let $\SimplexProd := \bigtimes_{\place \in \Set} \SimplexAt{\place}$
    be a \xnotion{polysimplex} in $\AffBuildProd$. For each place
    $\place\in\Set$, the stabilizer $\GroupAt{\place}$ of
    $\SimplexAt{\place}$ in $\UniOver{\FieldAt{\place}}$ is an
    open compact subgroup. For $\place\not\in\Set$, we put
    $\GroupAt{\place} := \UniOver{\OkaAt{\place}}$. Then
    $\AdeleGroup := \bigtimes_{\place} \GroupAt{\place}$ is an
    open subgroup of $\UniOver{\AdeleTo{\Set}}$.

    There is an obvious action of $\UniOver{\AdeleTo{\Set}}$ on
    $\AffBuildProd$: Components outside $\Set$ act trivially
    whereas a component corresponding to a place $\place\in\Set$ acts
    on the factor $\AffBuildAt{\place}$. Hence
    $
      \AffBuildProd = \UniOver{\AdeleTo{\Set}} \cdot \AffAppProd.
    $
    The stabilizer of $\SimplexProd$ is $\AdeleGroup$ whence
    the fiber of $\ProProd$ over $\SimplexProd$ is isomorphic to
    $
      \UniOver{\AdeleTo{\Set}} \rmod \AdeleGroup
    $
    which in turn is a discrete set since $\AdeleGroup$ is open.

    The group
    $\UniOver{\OkaSet} \leq \UniOver{\AdeleTo{\Set}}$ acts on $\UniOver{\AdeleTo{\Set}} \rmod \AdeleGroup$.
    Since $\OkaSet = \Field \cap \AdeleTo{\Set}$, Lemma~\ref{eins:kokompakt}
    implies that the double quotient
    $$
      \UniOver{\OkaSet} \lmod \UniOver{\AdeleTo{\Set}} \rmod \AdeleGroup
    $$
    is discrete and compact. Thus, it is finite.

    Therefore, the $\ProProd$--fiber over each polysimplex consists
    of finitely many $\UniOver{\OkaSet}$--orbits of cells.
    Now the claim is evident.
  \end{proof}
  \begin{Lemma}\label{Stabilisatoren:endlich}
    The group $\ChevOver{\OkaSet}$ acts on
    $\AffBuildProd$ with finite cell stabilizers.
  \end{Lemma}
  \begin{proof}
    The cell stabilizers of the action of $\ChevOver{\FieldAt{\place}}$ on $\AffBuildAt{\place}$
    are compact.
    Indeed, vertex stabilizers of this action are maximal compact subgroups of
    $\ChevOver{\FieldAt{\place}}$ \cite[Section~3.3]{Bruhat.Tits:1972}.
    Therefore, the stabilizer in $\ChevOver{\AdeleTo{\Set}}$ of a polysimplex in
    $\AffBuildProd$ is
    compact. The claim follows since $\ChevOver{\OkaSet}$ is a discrete subgroup
    of $\ChevOver{\AdeleTo{\Set}}$.
  \end{proof}

  Each standard apartment $\AffAppAt{\place}$ is a Euclidean space
  of dimension $\TorDim$. Within each of them, we choose a sector
  \begin{notationlist}
  \item
    $\SectorAt{\place}\subseteq\AffAppAt{\place}$ representing the
    fundamental chamber $\InfChamberAt{\place}$ of $\InfBuildAt{\place}$.
    We regard the cone point of $\SectorAt{\place}$ as the origin
    in $\AffAppAt{\place}$ turning the apartment into a Euclidean vector space.
    Moreover, we represent all roots in $\RootSystem$ as linear forms
    on $\AffAppAt{\place}$. Following the usual convention, we call
    those of them \notion{negative} that take negative values inside
    $\SectorAt{\place}$. Thus we are given a system
  \item
    $\NegativeRootsAt{\place}$ of negative roots in $\AffAppAt{\place}$.
    Considered as a subset of $\RootSystem$, it is independent of the place $\place$
    since all fundamental chambers $\InfChamberAt{\place}$ correspond
    to the same Borel subgroup scheme $\Borel$. Passing to a set of
    \xnotion{base roots},
    we obtain a system of coordinates
  \item
    \(
      \CoordinatesAt{\place} :=
        \bTupelOf{
           \AffRootAtInd{\place}{1},
           \ldots,
           \AffRootAtInd{\place}{\TorDim}
        }
    \)
    on $\AffAppAt{\place}$.
    With respect to these coordinates, the sector $\SectorAt{\place}$
    is given by
    $$
      \SectorAt{\place} =
      \bSetOf{ 
        \AffAppPointAt{\place} \in \AffAppAt{\place}
        \bsuchthatvrule
        \AffRootAtInd{\place}{\Index}( \AffAppPointAt{\place} ) \leq 0
        \;\;\forall \Index \in \SetOf{1,\ldots,\TorDim}
      }.
    $$
    Thus
    $$
      \AffAppProdTupel =
      \FamOf{
        \AffAppPointAt{\place}
      }{
        \place \in \Set
      }
      \mapsto
      \FamOf{
        \bTupelOf{
          \AffRootAtInd{\place}{1}(\AffAppPointAt{\place}),
          \ldots,
          \AffRootAtInd{\place}{\TorDim}(\AffAppPointAt{\place})
        }
      }{
        \place \in \Set
      }
    $$
    defines coordinates on $\AffAppProd$. Scaling the
    different $\place$--components appropriately, we arrange things
    such that the action of $\TorusOver{\OkaSet}$ leaves the map
  \item
    \ifNotationlisting
    \(
      \CoordinatesSum
      \mapcolon 
      \AffAppProdTupel =
      \FamOf{
        \AffAppPointAt{\place}
      }{
        \place \in \Set
      }
      \mapsto
      \bTupelOf{
        \sum_{\place \in \Set}
        \AffRootAtInd{\place}{1}(\AffAppPointAt{\place}),
        \ldots,
        \sum_{\place \in \Set}
        \AffRootAtInd{\place}{\TorDim}(\AffAppPointAt{\place})
      }
    \)
    \else
    \[
      \CoordinatesSum
      \mapcolon 
      \AffAppProdTupel =
      \FamOf{
        \AffAppPointAt{\place}
      }{
        \place \in \Set
      }
      \mapsto
      \bbTupelOf{
        \sum_{\place \in \Set}
        \AffRootAtInd{\place}{1}(\AffAppPointAt{\place}),
        \ldots,
        \sum_{\place \in \Set}
        \AffRootAtInd{\place}{\TorDim}(\AffAppPointAt{\place})
      }
    \]
    \fi
    invariant. We can do so because the idelic norm is identically $1$ on
    $\UnitsOf{\Field}$ by the
    product formula \cite[page~60]{Cassels.Froehlich:1967}.
    The coordinates $\AffRootAtInd{\place}{\Index}$
    correspond to base roots and, therefore, to characters
    $\CharacterInd{\Index} \mapcolon \Torus \rightarrow \Mult$. An element
    $\TorusElemAt{\place} \in \TorusOver{\FieldAt{\place}}$
    acts on $\AffAppAt{\place}$ by a translation whose
    $\Index$--coordinate is
    $\degat{\place}\place(\CharacterInd{\Index}(\TorusElemAt{\place}))$.

    We let $\BorelOver{\AdeleTo{\Set}}$ act on $\AffAppProd$ via the
    projection $\Borel \rightarrow \Torus$. This way, $\ProProd$
    becomes a $\BorelOver{\AdeleTo{\Set}}$--map and hence a
    $\BorelOver{\OkaSet}$--map. Since this action of $\BorelOver{\OkaSet}$
    on $\AffAppProd$
    factors through the
    torus $\Torus$ it leaves the map
    \item
      $\ProSum := \CoordinatesSum \circ \ProProd$ invariant.
  \end{notationlist}
  \begin{Lemma}\label{Operation:kokompakt}
    For any compact subset
    $\Kompaktum \subseteq \RRR^{\TorDim}$, the preimage
    $\ProSum^{-1}(\Kompaktum)$ contains a compact subset whose
    $\BorelOver{\OkaSet}$--translates cover $\ProSum^{-1}(\Kompaktum)$.
  \end{Lemma}
  \begin{proof}
    Dirichlet's Unit Theorem~\cite[page~72]{Cassels.Froehlich:1967}
    implies that $\TorusOver{\OkaSet}$ acts cocompactly on
    the kernel of $\CoordinatesSum$ whence it acts cocompactly
    on the preimage $\CoordinatesSum^{-1}(\Kompaktum)$ as well.
    So, let $\KompaktDomain \subseteq \AffAppProd$ be a compact set
    whose $\TorusOver{\OkaSet}$--translates
    cover $\CoordinatesSum^{-1}(\Kompaktum)$.

    By Lemma~\ref{Projektion:eigentlich}, we can find a compact subset
    $\AltKompaktDomain\subseteq\AffBuildProd$ whose $\UniOver{\OkaSet}$--translates
    cover $\ProProd^{-1}(\KompaktDomain)$. Then, the $\BorelOver{\OkaSet}$--translates
    of $\AltKompaktDomain$ cover $\ProSum^{-1}(\Kompaktum)$.
  \end{proof}

\section{Example: Rank--1--groups and trees}\label{Sec:SlZwei}
  It is as instructive as useful to treat the most simple case
  first: the Chevalley group $\SlInd{2}$. Serre gives a comprehensive discussion
  of this group and its associated building in \cite[II.1]{Serre:1980}.
  The group $\SDiagInd{2}$ of diagonal matrices with determinant~$1$
  is a maximal torus, and the group $\SUTriInd{2}$ of upper triangular matrices
  with determinant~$1$ is a Borel subgroup. Its unipotent part is the
  group $\UniTriInd{2}$ of strict upper triangular matrices all of whose diagonal
  entries equal $1$.
  The affine building $\AffBuildAt{\place}$ at the place~$\place$ is
  a regular tree of order $\CardOf{\fieldat{\place}}+1$: points in the
  link of a vertex correspond to points of the projective line over
  $\fieldat{\place}$.
  The standard apartment $\AffAppAt{\place}$ is a line. The projection
  map $\ProAt{\place}$ can be regarded as a \notion{height function}
  on $\AffBuildAt{\place}$ by identifying the apartment
  $\AffAppAt{\place}$ with the real line $\RRR$ via the negative base root
  $\AffRootAt{\place} \mapcolon \AffAppAt{\place} \rightarrow \RRR$.
  By scaling, as in the general case, we arrange that
  the action of $\SDiagIndOver{2}{\OkaSet}$ on the product
  $\AffBuildProd := \bigtimes_{\place \in \Set} \AffBuildAt{\place}$
  leaves the height
  function $\ProSum \mapcolon \AffBuildProd \rightarrow \RRR$ invariant.
  This situation has been discussed already in~\cite{Bux:1997}.
  Here we will treat it without making use of Bieri--Strebel theory.

  Trees are crucial for everything that follows. Therefore, we will
  repeatedly make use of the following lemma, which may look somewhat technical
  at a first glance. However, it describes a rather
  natural geometrical situation.
  \begin{Lemma}\label{BaumProdukt}
    Let
    \(
      \FamOf{
        \HightInd{\Index} \mapcolon \TreeInd{\Index} \rightarrow \RRR
      }{
        \Index\in\SetOf{1,\ldots,\MaxIndex}
      }
    \)
    be a family of locally
    finite simplicial trees $\TreeInd{\Index}$ with height
    functions $\HightInd{\Index}$. Suppose for every index
    $\Index$,
    \begin{enumerate}
      \item
        $\HightInd{\Index}$ maps the vertices of $\TreeInd{\Index}$ to
        a discrete subset of $\RRR$;
      \item
        there is exactly one descending end in $\TreeInd{\Index}$, ie, any two edge paths
        along which the height strictly decreases will eventually conincide; and
      \item
        each vertex in $\TreeInd{\Index}$ has degree $\geq 3$.
    \end{enumerate}
    So all descending paths eventually meet, and every vertex has a unique lower neighbor and at least two
    higher neighbors.

    Let $\TreeProd:=\TreeInd{1}\times\cdots\times\TreeInd{\MaxIndex}$ be the product of the trees $\TreeInd{\Index}$ and
    let $\HightSum \mapcolon \TreeProd \rightarrow \RRR$ be defined by    
    $$
      \HightSum \mapcolon \TreeProdTupel = 
      (\TreePointInd{1},\ldots,\TreePointInd{\MaxIndex})
      \mapsto
      \sum_{\Index = 1}^{\MaxIndex} \HightInd{\Index}( \TreePointInd{\Index} ).
    $$
    For every compact interval $\Intervall\subset\RRR$, put
    $\TreeProdSub{\Intervall} := \HightSum^{-1}(\Intervall)$.

    Then, for each compact interval $\Intervall$, the space
    $\TreeProdSub{\Intervall}$ is \notion{$(\MaxIndex-2)$--connected},
    ie, the homotopy groups $\FundamentalGroup{\Index}{\TreeProdSub{\Intervall}}$
    are trivial for $0\leq\Index\leq\MaxIndex-2$.

    Moreover, for any two intervals $\Intervall \subseteq \SupIntervall$,
    the inclusion
    $\TreeProdSub{\Intervall} \monorightarrow \TreeProdSub{\SupIntervall}$
    induces a non-trivial map
    $\RedHomTo{\MaxIndex-1}{\TreeProdSub{\Intervall}}
     \rightarrow
     \RedHomTo{\MaxIndex-1}{\TreeProdSub{\SupIntervall}}$
    in reduced homology.
  \end{Lemma}
  \begin{proof}
    The map $\HightSum$ is a \xnotion{Morse function} as defined
    in~\cite[Definition~2.2]{Bestvina.Brady:1997}. Its \xnotion{ascending} and
    \xnotion{descending} links in $\TreeProd$ are the joins
    of the ascending and descending links of the $\HightInd{\Index}$
    in the trees $\TreeInd{\Index}$, respectively. Thus, the descending
    links are points, and the ascending links are wedges of $(\MaxIndex-1)$--spheres.
    Hence ascending
    and descending links in $\TreeProd$ are $(\MaxIndex-2)$--connected.
    Then \cite[Corollary~2.6]{Bestvina.Brady:1997} implies that
    $\TreeProdSub{\Intervall}$ is $(\MaxIndex-2)$--connected
    for each interval $\Intervall$: The product $\TreeProd$
    could not be contractible otherwise.

    As for the second claim, recall that
    each tree $\TreeInd{\Index}$ has a unique descending end.
    Moving every point in $\TreeInd{\Index}$ with unit speed downhill toward
    this end defines a flow on $\TreeInd{\Index}$.
    We obtain a flow on $\TreeProd$ that moves all points
    in $\TreeProdSub{\Intervall}$ toward
    $\TreeProdSub{\SetOf{\min(\Intervall)}}$. This construction shows
    that $\TreeProdSub{\SetOf{\min(\Intervall)}}$ is a strong
    deformation retract of $\TreeProdSub{\Intervall}$.

    For $\Intervall\subseteq\SupIntervall$ the retraction
    $\TreeProdSub{\SupIntervall}\rightarrow\TreeProdSub{\SetOf{\min(\SupIntervall)}}$
    induces a map
    \(
      \TreeProdSub{\SetOf{\min(\Intervall)}}\rightarrow\TreeProdSub{\SetOf{\min(\SupIntervall)}}
    \)
    such that the following diagram
    \ifXYpic
    \[
      \xymatrix{
        {\RedHomTo{\MaxIndex-1}{\TreeProdSub{\Intervall}}}
        \ar[d] \ar@{=}[r]
      &
        {\RedHomTo{\MaxIndex-1}{\TreeProdSub{\SetOf{\min(\Intervall)}}}}
        \ar[d]
      \\
        {\RedHomTo{\MaxIndex-1}{\TreeProdSub{\SupIntervall}}}
        \ar@{=}[r]
      &
        {\RedHomTo{\MaxIndex-1}{\TreeProdSub{\SetOf{\min(\SupIntervall)}}}}
      \\
      }
    \]
    \else
    \[
      \begin{array}{ccc}
        \RedHomTo{\MaxIndex-1}{\TreeProdSub{\Intervall}}    & = & \RedHomTo{\MaxIndex-1}{\TreeProdSub{\SetOf{\min(\Intervall)}}} \\
        \downarrow                                          &   & \downarrow \\
        \RedHomTo{\MaxIndex-1}{\TreeProdSub{\SupIntervall}} & = & \RedHomTo{\MaxIndex-1}{\TreeProdSub{\SetOf{\min(\SupIntervall)}}} \\
      \end{array}
    \]
    \fi
    commutes.
    We will construct a sphere in $\TreeProdSub{\SetOf{\min(\Intervall)}}$
    which maps to a non-trivial embedded $(\MaxIndex-1)$--sphere in
    $\TreeProdSub{\SetOf{\min(\SupIntervall)}}$.
    This proves the claim because the latter sphere
    defines a non-trivial cycle that cannot be a boundary because
    there is no $\MaxIndex$--skeleton in
    $\TreeProdSub{\SetOf{\min(\SupIntervall)}}$.

    We choose a point
    $\TreeProdTupel = (\TreePointInd{1},\ldots,\TreePointInd{\MaxIndex})
    \in \TreeProd$
    with $\HightSum(\TreeProdTupel) < \min(\SupIntervall)$
    all of whose coordinates $\TreePointInd{\Index}\in\TreeInd{\Index}$
    are vertices.
    For each $\Index$, we choose two ascending rays
    $\LineIndDir{\Index}{+}$ and $\LineIndDir{\Index}{-}$ starting
    at $\TreePointInd{\Index}$ without a common initial segment---recall
    that every vertex has at least two higher neighbors.
    The union $\LineInd{\Index} := 
    \LineIndDir{\Index}{+} \cup \LineIndDir{\Index}{-}$ is a line
    in $\TreeInd{\Index}$. The distance of a point
    $\LinePointInd{\Index} \in \LineInd{\Index}$
    to the ``splitting vertex'' $\TreePointInd{\Index}$ is given by
    $\HightInd{\Index}(\LinePointInd{\Index}) - 
    \HightInd{\Index}(\TreePointInd{\Index})$. Hence
    the map $\LineProdTupel \mapsto
    \HightSum(\LineProdTupel) - \HightSum(\TreeProdTupel)$
    defines a norm on the product
    $\LineProd := \bigtimes_{\Index = 1}^{\MaxIndex} \LineInd{\Index}$.

    The sphere we wanted is the sphere of all
    points in $\LineProd$ whose norm is
    $\min(\Intervall) - \HightSum( \TreeProdTupel )$. The retraction
    shrinks it to the sphere of radius
    $\min(\SupIntervall) - \HightSum( \TreeProdTupel )$,
    which is still strictly positive.
  \end{proof}
  The rank--$1$--case is now easy since
  we can invoke K\,S Brown's celebrated criterion:
  \begin{cit}[{\cite[Remark~(2) to Theorem~2.2 and
              Theorem~3.2]{Brown:1987}}]\label{Kriterium:Brown}
    Let $\Group$ be a group, $\DirSet$ a directed set, and
    $\FamOf{\ComplexTo{\direl}}{\direl \in \DirSet}$ a
    directed system of \CW-complexes on which $\Group$ acts by
    cell permuting homeomorphisms such that the following hold:
    \begin{enumerate}
      \item\label{Kriterium:Brown:kokompakt}
        For each $\direl \in \DirSet$, the orbit space
        $\Group \lmod \ComplexTo{\direl}$ is compact.
      \item\label{Kriterium:Brown:Stabilisatoren}
        The stabilizer in $\Group$ of each $\celldim$--cell in $\ComplexTo{\direl}$
        is a group of type \F{\type - \celldim}.
      \item\label{Kriterium:Brown:Funktor}
        The continuous map 
        $\ComplexTo{\direla} \rightarrow \ComplexTo{\direlb}$
        indexed by $\direla \leq \direlb$ is $\Group$--equivariant.
      \item\label{Kriterium:Brown:Limes}
        The limit of the directed system of homotopy groups
        $\FamOf{\FundamentalGroup{\Index}{\ComplexTo{\direl}}}{\direl \in \DirSet}$
        vanishes for $\Index < \type$.
    \end{enumerate}
    Then, $\Group$ is of type \FP{\type} if and only if for all
    $\Index<\type$, the directed system of reduced homology groups
    $\sFamOf{\RedHomTo{\Index}{\ComplexTo{\direl}}}{\direl \in \DirSet}$
    is \notion{essentially trivial}, ie, for each $\direla \in \DirSet$,
    there is $\direlb \geq \direla$ such that the map
    $\RedHomTo{\Index}{\ComplexTo{\direla}} \rightarrow
    \RedHomTo{\Index}{\ComplexTo{\direlb}}$ induced by
    $\ComplexTo{\direla} \rightarrow \ComplexTo{\direlb}$
    is trivial.

    Moreover, $\Group$ is of type \F{\type} if and only if the directed
    system
    of homotopy groups $\FamOf{\FundamentalGroup{\Index}{\ComplexTo{\direl}}}{\direl \in \DirSet}$ is essentially trivial
    for all $\Index<\type$.
  \end{cit}
  \begin{cor}
    \label{Brown:Simple_Criterion}
    Let $\Group$ act cocompactly by cell-permuting homeomorphisms on an $(\type-1)$--connected
    CW-complex $\Complex$ such that the stabilizer of each cell is finite. Then $\Group$ is
    of type \F{\type}.
  \end{cor}
  \begin{proof}
    Take a directed set consisting of just one element, assign $\Complex$
    as the corresponding complex, and observe that the identity map induces
    trivial maps in homotopy groups in those dimensions where these groups vanish.
    Since finite groups are of type \F{\infty}, the claim follows from
    Brown's Criterion.
  \end{proof}
  \begin{rem}\label{Condition_Four}
    In our applications, the direct limit of the spaces $\ComplexTo{\direl}$
    will be the union of these spaces. Usually, it will be contractible which
    then implies that
    the limit of the directed system of homotopy groups
    $\FamOf{\FundamentalGroup{\Index}{\ComplexTo{\direl}}}{\direl \in \DirSet}$
    vanishes for $\Index < \type$.
  \end{rem}
  \begin{cor}\label{SlZwei}
    The group $\SUTriIndOver{2}{\OkaSet}$ is of type
    \F{\CardOf{\Set}-1} but not of type \FP{\CardOf{\Set}}.
  \end{cor}
  \begin{proof}
    We apply Brown's Criterion. The set of all compact intervals
    in $\RRR$ is a directed set ordered by inclusion, and we are
    looking for a family of cocompact $\SUTriIndOver{2}{\OkaSet}$--\CW--complexes
    over this system. $\SUTriIndOver{2}{\OkaSet}$ acts on the
    product of trees $\AffBuildProd$ with $\ProSum$ as an invariant
    height function.
    Hence for each compact interval $\Intervall$, the preimage
    $\Sub{\AffBuildProd}{\Intervall}:=\ProSum^{-1}(\Intervall)$ is
    a $\SUTriIndOver{2}{\OkaSet}$--complex. This defines our
    directed system with inclusions as
    continuous, $\SUTriIndOver{2}{\OkaSet}$--equivariant maps.

    The hypotheses of Brown's Criterion are satisfied:
    The action is by cell-permuting homeomorphisms, it is cocompact by
    Lemma~\ref{Operation:kokompakt}, 
    cell stabilizers are even finite by
    Lemma~\ref{Stabilisatoren:endlich}, and
    condition~\ref{Kriterium:Brown:Limes} is satisfied because the
    limit $\AffBuildProd$ of all $\Sub{\AffBuildProd}{\Intervall}$
    is contractible.
    
    The height function $\ProSum$ can be regarded as a
    sum of height function defined on the factors $\AffBuildAt{\place}$ such
    that we are in the setting of Lemma~\ref{BaumProdukt}: the descending
    end in $\AffBuildAt{\place}$ is the unique chamber at infinity (in this
    case just a point) stabilized by $\SUTriIndOver{2}{\FieldAt{\place}}$.
    This completes the proof.
  \end{proof}

  What can we say about the related group scheme $\GlInd{2}$? The short
  exact sequence
  $\SlInd{2} \monorightarrow \GlInd{2} \epirightarrow \Mult$
  with the determinant as the projection map induces by restriction
  a short exact sequence
  \ifXYpic
  \[
    \xymatrix{
      {\SUTriIndOver{2}{\OkaSet}}
      \ar@{ >->}[r]
    &
      {\UTriIndOver{2}{\OkaSet}}
      \ar@{->>}[r]
    &
      {\UnitsOf{\OkaSet}}
    }
  \]
  \else
  \[
    \SUTriIndOver{2}{\OkaSet} \monorightarrow
    \UTriIndOver{2}{\OkaSet} \epirightarrow
    \UnitsOf{\OkaSet}
  \]
  \fi
  whence $\UTriIndOver{2}{\OkaSet}$ inherits all finiteness properties
  of $\SUTriIndOver{2}{\OkaSet}$ since $\UnitsOf{\OkaSet}$ is of type
  \F{\infty}. However, $\UTriIndOver{2}{\OkaSet}$ might even exhibit stronger
  finiteness properties, but we can rule out this possibility:
  \begin{rem}\label{GlZwei}
    The group $\UTriIndOver{2}{\OkaSet}$ is of type
    \F{\CardOf{\Set}-1} but not of type \FP{\CardOf{\Set}}.
  \end{rem}
  \begin{proof}
    Passing to projective groups, we obtain the following commutative
    diagram all of whose rows and columns are short exact sequences of
    groups:
    \ifXYpic
    \[
      \xymatrix{
        {\SetOf{-1,1}}
        \ar@{ >->}[r] \ar@{ >->}[d]_{{\cdot \III_2}}
      &
        {\UnitsOf{\OkaSet}}
        \ar@{->>}[r]^{(\cdot)^2} \ar@{ >->}[d]_{{\cdot \III_2}}
      &
        {\UnitSquaresOf{\OkaSet}}
        \hbox to 0pt{\(
          \null:=\SetOf{ \function^2 \suchthatvrule \function \in \UnitsOf{\OkaSet}}
        \)}
        \ar@{ >->}[d]
      \\
        {\SUTriIndOver{2}{\OkaSet}}
        \ar@{ >->}[r]  \ar@{->>}[d]
      &
        {\UTriIndOver{2}{\OkaSet}}
        \ar@{->>}[r]  \ar@{->>}[d]
      &
        {\UnitsOf{\OkaSet}}
        \ar@{->>}[d]
      \\
        {\PSUTriIndOver{2}{\OkaSet}}
        \ar@{ >->}[r]
      &
        {\PUTriIndOver{2}{\OkaSet}}
        \ar@{->>}[r]
      & 
        {\UnitsOf{\OkaSet}\rmod\UnitSquaresOf{\OkaSet}}
      \\
      }
    \]
    \else
    \[
      \begin{array}{ccccc}
        \SetOf{-1,1}               & \monorightarrow & \UnitsOf{\OkaSet}         & \stackrel{(\cdot)^2}{\epirightarrow} & 
                                         \UnitSquaresOf{\OkaSet} \vcenter{\rlap{$:=\SetOf{ \function^2 \suchthatvrule \function \in \UnitsOf{\OkaSet}}$}}   \\
        \ldec{\cdot \III_2}        &                 & \ldec{\cdot \III_2}       &                & \ldec{}                                       \\
        \SUTriIndOver{2}{\OkaSet}  & \monorightarrow & \UTriIndOver{2}{\OkaSet}  & \epirightarrow & \UnitsOf{\OkaSet}                             \\
        \ldec{}                    &                 & \ldec{}                   &                & \ldec{}                                       \\
        \PSUTriIndOver{2}{\OkaSet} & \monorightarrow & \PUTriIndOver{2}{\OkaSet} & \epirightarrow & \UnitsOf{\OkaSet}\rmod\UnitSquaresOf{\OkaSet} \\
      \end{array}
    \]
    \fi
    Consider the bottom row first. The factor
    $\UnitsOf{\OkaSet}\rmod\UnitSquaresOf{\OkaSet}$
    on the right is an abelian
    torsion group which is finitely generated by Dirichlet's Unit Theorem.
    Hence $\UnitsOf{\OkaSet}\rmod\UnitSquaresOf{\OkaSet}$
    is finite. Therefore, the other two groups in this row
    enjoy the same finiteness properties.
    Then, this also holds for their extensions in the
    middle row because the kernels on top are of type \F{\infty}.
  \end{proof}
  \begin{rem}\label{GlDrei}
    Of course, the same argument implies that
    $\SUTriIndOver{\Rank}{\OkaSet}$ and $\UTriIndOver{\Rank}{\OkaSet}$
    enjoy identical finiteness properties.
  \end{rem}

\section{Higher ranks---an algebraic prelude}\label{BeginProof}
  Perhaps the most striking consequence of Theorem~A is
  that the finiteness length of $\BorelOver{\OkaSet}$ does not depend
  of the rank of the Chevalley group $\Chev$. In this section, 
  a simple algebraic explanation for the group scheme $\SUTriInd{\Rank}
  \subseteq \SlInd{\Rank}$ is given.
  \begin{prop}\label{GruppenRetrakt}
    Suppose
    \(
      \RetractDiagram{\Group}{\FaktorGroup}
    \)
    is a \notion{retract diagram} of groups, ie, the composition of
    arrows is the identity on $\FaktorGroup$. Then $\FaktorGroup$ inherits
    all finiteness properties of $\Group$.
  \end{prop}
  \begin{proof}
    Finite generation is trivial, finite presentability is easy and
    dealt with in~\cite[Lemma~1.3]{Wall:1965} where it is attributed to
    J\,R Stallings. Thus, it suffices to treat the homological finiteness
    properties starting with \FP{2}. It would be possible to cite
    \cite[($\zeta$), page~280]{Aberg:1986}, but \AA{}berg is merely hinting
    at the argument. The following claims to be what he had in mind.

    The key observation is, that functors and cofunctors both preserve retract
    diagrams. So for each index set $\IndexSet$, consider the functor
    that assigns to a group $\Group$ the pair
    $\TupelOf{ \Group, \bigtimes_\IndexSet \ZZZ \Group }$
    where we regard $\bigtimes_\IndexSet \ZZZ \Group$ as a
    $\ZZZ \Group$--module. Applying the homology functor
    $\HomTo{\celldim}{-,-}$ yields the following retract diagram:
    \[
      \RetractDiagram{
        \HomTo{\celldim}{\Group, \bigtimes_\IndexSet \ZZZ \Group}
      }{
        \HomTo{\celldim}{\FaktorGroup, \bigtimes_\IndexSet \ZZZ \FaktorGroup}
      }
    \]
    Hence
    $\HomTo{\celldim}{\FaktorGroup,\bigtimes_\IndexSet\ZZZ\FaktorGroup}$
    vanishes whenever
    $\HomTo{\celldim}{\Group,\bigtimes_\IndexSet\ZZZ\Group}$
    does.

    The claim now follows by means of the Bieri--Eckmann Criterion:
    A finitely generated group $\FaktorGroup$ is of type \FP{\type}
    if and only if
    $\HomTo{\celldim}{\FaktorGroup, \bigtimes_\IndexSet\ZZZ\FaktorGroup} = 0$
    for all index sets $\IndexSet$ and all
    $\celldim \in \SetOf{\One,\ldots,\type-\One}$
    \cite[Proposition~1.2 and the equation above
                     Theorem~2.3]{Bieri.Eckmann:1974}.
  \end{proof}
  A different proof, based not on the Bieri--Eckmann Criterion but
  on Brown's Criterion, can be found in \cite[Remark 3.3]{Bux:2002}.
  \begin{cor}
    The groups $\UTriIndOver{\Rank}{\OkaSet}$ and
    $\SUTriIndOver{\Rank}{\OkaSet}$ are not of type \FP{\CardOf{\Set}}.
  \end{cor}
  \begin{proof}
    We confine ourselves to $\Rank=3$. The group $\UTriInd{2}$
    embeds into $\UTriInd{3}$ like this:
    \[
      \UTriInd{2}
      \isom
      \left(\begin{array}{ccc}
          * & * & 0 \\
          0 & * & 0 \\
          0 & 0 & 1 \\
      \end{array}\right)
      \subseteq
      \UTriInd{3}
      .
    \]
    This way, we recognize $\UTriInd{2}$ as a retract of $\UTriInd{3}$.
    Hence, the preceding Proposition~\ref{GruppenRetrakt} applies and
    the claim follows from
    Corollary~\ref{SlZwei} and
    the Remarks~\ref{GlZwei}
    and~\ref{GlDrei}.
  \end{proof}

\section{A geometric version of the argument}
  \label{GeometrischeVersion}
  \begin{theorem}\label{ObereAbschaetzung}
    The group $\BorelOver{\OkaSet}$ is not of type \FP{\CardOf{\Set}}.
  \end{theorem}
  The entire section is devoted to the proof of this theorem.
  The reasoning can be viewed as
  a geometric interpretation of the argument presented in the preceding
  section. Let us start with a brief outline: Each of the affine buildings
  $\AffBuildAt{\place}$ contains a tree $\TreeAt{\place}$ as a
  retract. Hence the product $\AffBuildProd$ contains a product
  $\TreeProd$ of trees as a retract. We will find a directed system of
  subspaces in $\AffBuildProd$ satisfying the hypotheses of Brown's
  Criterion. So we only have to prove that the induced system of reduced
  homology groups is not essentially trivial. Finally, using the retraction map,
  we pass to a corresponding system of subspaces in $\TreeProd$ where
  Lemma~\ref{BaumProdukt} applies.

  Let us call an affine apartment in $\AffBuildAt{\place}$ a \notion{layer}
  if it contains the fundamental chamber at infinity
  $\InfChamberAt{\place}$.
  The base root
  $\AffRootAt{\place} := \AffRootAtInd{\place}{1}$ defines
  half apartments in $\AffBuildAt{\place}$ by
  $$
    \AffAppAtOf{\place}{\real}
    :=
    \bSetOf{
      \AffAppPointAt{\place} \in \AffAppAt{\place}
      \bsuchthatvrule
      \AffRootAt{\place}( \AffAppPointAt{\place} ) \leq \real
    }
    .
  $$
  Call an apartment in $\AffBuildAt{\place}$ \notion{special}
  if it contains such a half apartment. To put it in a slightly different way:
  The base root defines a half apartment at infinity
  $\InfRootAt{\place} \subseteq \InfBuildAt{\place}$,
  which contains the fundamental chamber
  $\InfChamberAt{\place}$. An affine apartment in $\AffBuildAt{\place}$
  is special if and only if it contains $\InfRootAt{\place}$.
  Obviously, every special apartment is a layer.

  The map
  $\AffHightAt{\place} := \AffRootAt{\place} \circ \ProAt{\place}$
  restricts to an affine map
  on every layer---hence on every special apartment. 
  Two special apartments
  $\AffAppAtInd{\place}{1}$ and $\AffAppAtInd{\place}{2}$
  intersect in a convex set, which contains a subset of the
  form $\AffAppAtOf{\place}{\lowreal}$.
  Hence,
  $$
    \AffAppAtInd{\place}{1} \cap \AffAppAtInd{\place}{2}
    =
    \SetOf{
      \AffAppPoint \in \AffAppAtInd{\place}{1}
      \suchthatvrule
      \AffHightAt{\place}( \AffAppPoint ) \leq \real
    }
    =
    \SetOf{
      \AffAppPoint \in \AffAppAtInd{\place}{2}
      \suchthatvrule
      \AffHightAt{\place}( \AffAppPoint ) \leq \real
    }
  $$
  where
  $
    \real
    =
    \max \AffHightAt{\place}(
      \AffAppAtInd{\place}{1} \cap \AffAppAtInd{\place}{2}
    )
  $.
  Thus we conclude:
  \begin{observation}\label{SpezielleApartments}
    The union of all special apartments in $\AffBuildAt{\place}$ is
    a subcomplex isometric to a product
    $\TreeAt{\place} \times \RRR^{\TorDim - 1}$ where
    $\TreeAt{\place}$ is a tree. The projection onto the second
    factor $\RRR^{\TorDim-1}$ is defined by
    $
      \bTupelOf{
        \AffRootAtInd{\place}{2} \circ \ProAt{\place}
        ,\ldots,
        \AffRootAtInd{\place}{\TorDim} \circ \ProAt{\place}
      }
    $.
    In particular, the fiber over each tuple
    $\TupelOf{ \RealInd{2}, \ldots ,\RealInd{\TorDim} }$
    is a tree on which $\AffHightAt{\place}$ induces a height function
    $\TreeHightAt{\place} \mapcolon \TreeAt{\place} \rightarrow \RRR$.\qed
  \end{observation}
    In~\cite[Chapter~10.2]{Ronan:1989}, M Ronan gives two constructions
    for the tree $\TreeAt{\place}$. The equivalence of these two constructions
    underlies the following argument.

  \begin{Lemma}\label{SpeziellesApartmentExistiert}
    For each layer $\AffAppAtInd{\place}{'}$, there is a special
    apartment $\AffAppAtInd{\place}{\mathrm{s}}$ such that
    $
      \AffHightAt{\place}(
        \AffAppAtInd{\place}{'} \cap \AffAppAtInd{\place}{\mathrm{s}}
      )
    $
    is unbounded.
  \end{Lemma}
  \begin{proof}
    We argue within the spherical building $\InfBuildAt{\place}$.
    The apartment at infinity $\InfAppAtInd{\place}{'}$ corresponding to
    $\AffAppAtInd{\place}{'}$ contains $\InfChamberAt{\place}$.
    The root $\InfRootAt{\place}$ determines a codimension $1$ face
    of the fundamental chamber $\InfChamberAt{\place}$. Let $\InfChamberAtInd{\place}{'}$ be
    the other neighboring chamber in $\InfAppAtInd{\place}{'}$. Then there is
    a unique apartment containing the chamber $\InfChamberAtInd{\place}{'}$
    and the half apartment $\InfRootAt{\place}$. The corresponding affine apartment
    $\AffAppAtInd{\place}{\mathrm{s}}$
    satisfies our needs since $\AffHightAt{\place}$ is unbounded on any
    sector representing $\InfChamberAtInd{\place}{'}$ in
    $\AffAppAtInd{\place}{\mathrm{s}}$.
  \end{proof}
  \begin{Lemma}\label{ProjektionWohldefiniert}
    Let $\AffAppAtInd{\place}{1}$ and $\AffAppAtInd{\place}{2}$
    be special apartments and $\AffAppAtInd{\place}{'}$ be a layer.
    Moreover, let $\real$ be a real number and let
    $
      \AffAppPointAtInd{\place}{\Index}
      \in
      \AffAppAtInd{\place}{\Index} \cap \AffAppAtInd{\place}{'}
    $
    be two points such that
    $
      \real = \AffHightAt{\place}( \AffAppPointAtInd{\place}{\Index} )
    $.
    Then
    $
      \real \in \AffHightAt{\place}(
        \AffAppAtInd{\place}{1} \cap \AffAppAtInd{\place}{2}
        )
    $.
  \end{Lemma}
  \begin{proof}
    Both intersections
    $\AffAppAtInd{\place}{\Index} \cap \AffAppAtInd{\place}{'}$
    are convex and contain a sector representing
    $\InfChamberAt{\place}$. Let $\SectorAtInd{\place}{\Index}$ be the
    corresponding parallel sector with cone point~$\AffAppPointAtInd{\place}{\Index}$.
    Both sectors $\SectorAtInd{\place}{1}$ and $\SectorAtInd{\place}{2}$
    are simplicial cones with a codimension~$1$~face restricted
    to which $\AffHightAt{\place}$ is identically $\real$. Hence, these two faces
    intersect within $\AffAppAtInd{\place}{'}$ and their intersection is
    a simplicial cone of codimension~$1$, which is contained in
    $\AffAppAtInd{\place}{1} \cap \AffAppAtInd{\place}{2}$.
  \end{proof}
  With the aid of the two lemmas above, we can see a projection
  $\RetraktionAt{\place} \mapcolon \AffBuildAt{\place} \rightarrow \TreeAt{\place}$.
  \begin{Lemma}\label{BaumRetrakt}
    There is a continuous projection map
    \(
      \RetraktionAt{\place}
      \mapcolon \AffBuildAt{\place} \rightarrow \TreeAt{\place}
    \)
    compatible with $\AffHightAt{\place}$, ie, the diagram
    \ifXYpic
    \[
      \xymatrix{
        {\AffBuildAt{\place}}
        \ar[r] \ar[d]_{\AffHightAt{\place}}
      &
        {\TreeAt{\place}}
        \ar[d]^{\TreeHightAt{\place}}
      \\
        {\RRR}
        \ar@{=}[r]
      &
        {\RRR}
      \\
      }
    \]
    \else
    \[
      \begin{array}{ccc}
        \AffBuildAt{\place} & \rightarrow & \TreeAt{\place}\\
        \ldec{\AffHightAt{\place}} &   & \ldec{\TreeHightAt{\place}}\\
        \RRR                & =           & \RRR
      \end{array}
    \]
    \fi
    commutes where $\TreeHightAt{\place}$ is as in
    Observation~\ref{SpezielleApartments}.
  \end{Lemma}
  \begin{proof}
    There is a $1$--$1$--correspondence
    $$
      \SetOf{
        \text{special apartments in\ }
        \AffBuildAt{\place}
      }
      \leftrightarrow
      \SetOf{
        \begin{array}{l}
          \text{lines in\ }
          \TreeAt{\place}
          \text{\ that\ }
          \TreeHightAt{\place}\\
          \text{maps isometrically to\ }
          \RRR
        \end{array}
      }
      =:
      \TheSetOfLines
      .
    $$
    According to Lemma~\ref{SpeziellesApartmentExistiert}, for each layer,
    there is a line onto which the layer can be projected in a way compatible
    with $\AffHightAt{\place}$. These projection maps agree where layers
    intersect by Lemma~\ref{ProjektionWohldefiniert}. Hence we have
    defined a projection map on $\AffBuildAt{\place}$ since the
    affine building is the union of all layers.

    This projection map is continuous when restricted to a layer since
    the tree $\TreeAt{\place}$ carries the weak topology with respect to
    the lines in $\TheSetOfLines$. Thus, the projection map is continuous because the building
    $\AffBuildAt{\place}$ carries the weak topology with respect to the layers.
  \end{proof}
  To see that this projection is a retraction, we need to find a
  continuous section. However, these exist in abundance by
  Observation~\ref{SpezielleApartments}.

  Consider the product
  \begin{notationlist}
    \item
      $\TreeProd := \bigtimes_{\place \in \Set} \TreeAt{\place}$,
      on which the height function
    \item
      $\TreeHightSum \mapcolon \TreeProd \rightarrow \RRR$
      is defined by
      $$
        \TreeProdTupel =
        \FamOf{ \TreePointAt{\place} }{\place \in \Set}
        \mapsto
        \sum_{\place \in \Set} \TreeHightAt{\place}(
          \TreePointAt{\place}
        ).
      $$
      There is the projection map
    \item
      $\RetraktionProd \mapcolon \AffBuildProd \rightarrow \TreeProd$
      that admits a lot of sections, which are parameterized by tuples
      $
        \TupelOf{
          \RealAtInd{\place}{\Index}
          \suchthatvrule 
          \place \in \Set ,\;
          \Index \in \SetOf{ 2,\ldots,\TorDim }
        }
      $
      of real numbers.
      Hence, $\TreeProd$ is a retract of $\AffBuildProd$.
  \end{notationlist}
  \begin{observation}\label{Beob:Kommutativ}
    The following diagram commutes:
    \ifXYpic
    \[
      \xymatrix{
        {\AffBuildProd}
        \ar@{->>}@<0.75mm>[r]  \ar[d]_{\ProSum}
      &
        {\TreeProd}
        \ar@{ >->}@<0.75mm>[l]  \ar[d]^{\TreeHightSum}
      \\
        {\RRR^{\TorDim}}
        \ar@{->>}[r]
      &
        {\RRR}
      }
    \]
    \else
    \[
      \begin{array}{ccc}
        \AffBuildProd  & \retraktarrows & \TreeProd \\
        \ldec{\ProSum} &                  & \ldec{\TreeHightSum}\\
        \RRR^{\TorDim} & \rightarrow      & \RRR
      \end{array}
    \]
    \fi
    Here, the arrow in the bottom row is the projection onto the first
    coordinate.\qed
  \end{observation}

  Finally, we can turn to finiteness properties of $\BorelOver{\OkaSet}$.
  We want to apply Brown's Criterion; thus, we specify a
  directed system of cocompact $\BorelOver{\OkaSet}$--subcomplexes in $\AffBuildProd$:
  A \notion{brick} is a product
  \(
    \IntervallProd :=
    \IntervallInd{1}
    \times \cdots \times
    \IntervallInd{\TorDim}
    \subset \RRR^{\TorDim}
  \)
  of compact intervals. The set of bricks is
  a directed set ordered by inclusion. By Lemma~\ref{Operation:kokompakt}, the family
  $\Sub{\AffBuildProd}{\IntervallProd} := \ProSum^{-1}( \IntervallProd )$
  of preimages of bricks
  is a directed system of cocompact
  $\BorelOver{\OkaSet}$--complexes.
  \begin{Lemma}\label{Lemma:Non_Trivial}
    The system~$\RedHomTo{\CardOf{\Set}-1}{
      \Sub{\AffBuildProd}{\IntervallProd}
    }$
    of reduced homology groups is not essentially trivial.
  \end{Lemma}
  \begin{proof}
    Given a brick $\IntervallProd$, let
    \(
      \SupIntervallProd :=
      \SupIntervallInd{1}
      \times \cdots \times
      \SupIntervallInd{\TorDim}
    \)
    be a brick containing $\IntervallProd$, ie,
    $\IntervallInd{\Index} \subseteq \SupIntervallInd{\Index}$
    for all $\Index$.
    Choose a tuple
    \(
      \TupelOf{
        \RealAtInd{\place}{\Index}
        \suchthatvrule 
        \place \in \Set ,\;
        \Index \in \SetOf{ 2,\ldots,\TorDim }
      }
    \)
    of real numbers such that
    \(
      \sum_{\place \in \Set} \RealAtInd{\place}{\Index}
      \in \IntervallInd{\Index}
    \)
    for $2 \leq \Index \leq \TorDim$.
    This defines a section of $\RetraktionProd$, which restricts to a section
    \[
      \Sub{\TreeProd}{\IntervallInd{1}} :=
      \TreeHightSum^{-1}( \IntervallInd{1} )
      \rightarrow
      \bbSub{\AffBuildProd}{
        \IntervallInd{1}\times
        \bbSetOf{\sum_{\place \in \Set} \RealAtInd{\place}{2}}
        \times\cdots\times
        \bbSetOf{\sum_{\place \in \Set} \RealAtInd{\place}{\TorDim}}
      }
      \subseteq
      \Sub{\AffBuildProd}{\IntervallProd}
      .
    \]
    Thus, Observation~\ref{Beob:Kommutativ} implies that the preimage
    \(
      \Sub{\TreeProd}{\IntervallInd{1}}
    \)
    is a retract of
    $\Sub{\AffBuildProd}{\IntervallProd}$
    and that the diagram
    \ifXYpic
    \[
      \xymatrix{
        {\Sub{\AffBuildProd}{\IntervallProd}}
        \ar@{->>}@<0.75mm>[r]  \ar[d]
      &
        {\Sub{\TreeProd}{\IntervallInd{1}}}
        \ar@{ >->}@<0.75mm>[l]  \ar[d]
      \\
        {\Sub{\AffBuildProd}{\SupIntervallProd}}
        \ar@{->>}@<0.75mm>[r]
      &
        {\Sub{\TreeProd}{\SupIntervallInd{1}}}
        \ar@{ >->}@<0.75mm>[l]
      }
    \]
    \else
    \[
      \begin{array}{ccc}
        \Sub{\AffBuildProd}{\IntervallProd} & 
          \retraktarrows &
            \Sub{\TreeProd}{\IntervallInd{1}}\\
        \downarrow & & \downarrow \\
        \Sub{\AffBuildProd}{\SupIntervallProd} &
          \retraktarrows &
            \Sub{\TreeProd}{\SupIntervallInd{1}}
      \end{array}
    \]
    \fi
    commutes. Passing to homology, we obtain the following
    commutative diagram:
    \ifXYpic
    \[
      \xymatrix{
        {\RedHomTo{\CardOf{\Set}-1}{\Sub{\AffBuildProd}{\IntervallProd}}}
        \ar@{->>}@<0.75mm>[r]  \ar[d]
      &
        {\RedHomTo{\CardOf{\Set}-1}{\Sub{\TreeProd}{\IntervallInd{1}}}}
        \ar@{ >->}@<0.75mm>[l]  \ar[d]
      \\
        {\RedHomTo{\CardOf{\Set}-1}{\Sub{\AffBuildProd}{\SupIntervallProd}}}
        \ar@{->>}@<0.75mm>[r]
      &
        {\RedHomTo{\CardOf{\Set}-1}{\Sub{\TreeProd}{\SupIntervallInd{1}}}}
        \ar@{ >->}@<0.75mm>[l]
      }
    \]
    \else
    \[
      \begin{array}{ccc}
        \RedHomTo{\CardOf{\Set}-1}{\Sub{\AffBuildProd}{\IntervallProd}} & 
          \retraktarrows & 
            \RedHomTo{\CardOf{\Set}-1}{\Sub{\TreeProd}{\IntervallInd{1}}} \\
        \downarrow & & \downarrow\\
        \RedHomTo{\CardOf{\Set}-1}{\Sub{\AffBuildProd}{\SupIntervallProd}}&
          \retraktarrows &
            \RedHomTo{\CardOf{\Set}-1}{\Sub{\TreeProd}{\SupIntervallInd{1}}}
      \end{array}
    \]
    \fi
    The right vertical arrow is non-trivial by Lemma~\ref{BaumProdukt}
    whence the left vertical arrow cannot be trivial either since the right hand side is
    a retract of the left hand side.
  \end{proof}
  To finish the proof of Theorem~\ref{ObereAbschaetzung}, observe that the preimages
  $\Sub{\AffBuildProd}{\IntervallProd}$ exhaust $\AffBuildProd$, which is contractible.
  The action of $\BorelOver{\OkaSet}$ on $\AffBuildProd$ is by cell-permuting homeomorphisms,
  and cell-stabilizers are finite by Lemma~\ref{Stabilisatoren:endlich}.
  Thus, in view of Brown's Criterion~\ref{Kriterium:Brown},
  Lemma~\ref{Lemma:Non_Trivial} completes the proof of Theorem~\ref{ObereAbschaetzung}.\qed

\section{The Moufang property}\label{Moufang}
    Fix a building $\Build$ and an apartment $\Apartment$ therein.
    For any half apartment $\roota$, let $-\roota$ denote the
    complementary half apartment.
    Two half apartments $\roota$ and $\rootb$ in $\Apartment$ are
    \notion{prenilpotent} if 
    $\roota \cap \rootb$ contains a chamber and
    $-\roota \cap -\rootb$ contains a chamber, too.
    In this case, we set
    \ifNotationlisting
    \begin{notationlist}
      \item
        $\AllBetween{\roota}{\rootb} :=
        \SetOf{\rootc \suchthatvrule \roota\cap\rootb \subseteq \rootc, 
                       -\roota\cap-\rootb \subseteq -\rootc}
        $ and
      \item
        $\InBetween{\roota}{\rootb} := \AllBetween{\roota}{\rootb}
                                   \setminus
                                   \SetOf{\roota,\rootb}.
        $
    \end{notationlist}
    \else
      \[
        \AllBetween{\roota}{\rootb} :=
        \SetOf{\rootc \suchthatvrule \roota\cap\rootb \subseteq \rootc, 
                       -\roota\cap-\rootb \subseteq -\rootc}
       \]
       and
       \[
         \InBetween{\roota}{\rootb} := \AllBetween{\roota}{\rootb}
                                   \setminus
                                   \SetOf{\roota,\rootb}.
       \]
    \fi
    $\Build$ is called \notion{Moufang} if one can associate
    \notion{root groups}
    \begin{notationlist}
      \item
        $\RootGroupAt{\root}$ of automorphisms of the
        building $\Build$ to the half apartments $\root \subseteq \Apartment$
        such that the following axioms hold:
    \end{notationlist}
    \begin{propositionsystem}{M}
      \item\label{MoufangEins}
        $\RootGroupAt{\root}$ fixes every chamber in $\root$, and for each
        panel $\panel\in\Boundary{\root}$ (a \notion{panel} is a codimension $1$ face of a chamber),
        $\RootGroupAt{\root}$ acts simply-transitively on
        the set of chambers in $\Star{\panel}$ but outside $\root$.
      \item\label{MoufangZwei}
        For each prenilpotent pair $\SetOf{\roota,\rootb}$, we have
        $[ \RootGroupAt{\roota}, \RootGroupAt{\rootb} ]
        \leq \RootGroupAt{\InBetween{\alpha}{\beta}}$.
        Here,
        $\RootGroupAt{\InBetween{\roota}{\rootb}}$ denotes the group
        generated by all $\RootGroupAt{\rootc}$ with
        $\rootc \in \InBetween{\roota}{\rootb}$.
      \item\label{MoufangDrei}
        For each $\RootGroupEl \in \RootGroupAt{\root} \setminus \{1\}$,
        there is an element $\mouf{\RootGroupEl} \in \RootGroupAt{-\root}
        \RootGroupEl \RootGroupAt{-\root}$ stabilizing $\Apartment$.
      \item\label{MoufangVier}
        For $\MoufEl = \mouf{u}$ as in~\ref{MoufangDrei}, we have
        $\MoufEl \RootGroupAt{\rootb} \MoufEl^{-1} = 
        \RootGroupAt{\MoufEl\rootb}$ for all half apartments $\rootb$.
    \end{propositionsystem}
  \begin{fact}
    The affine building associated to a Chevalley group $\Chev$
    over a local function field is Moufang.
    In particular, this holds for the buildings
    $\AffBuildAt{\place}$.
  \end{fact}
  This is ``well-known'' to those who, well, know it. However, there seems
  to be no explicit reference for this fact in the literature. For this reason,
  an outline of the argument is included.
  \begin{sketch}
    By Hensel's Lemma or
    \cite[Theorem~8, page~20]{Weil:1973}, any
    local function field is isomorphic, as a field with a valuation,
    to a field of Laurent series over a finite field. This can be regarded as
    the completion of the field of rational functions over the same
    finite field. The associated buildings are isomorphic. Thus,
    we may confine ourselves to the case of a rational function field.
    As P Abramenko observes in \cite[page~19]{Abramenko:1996},
    the corresponding building is isomorphic to the positive partner
    within the twin building of $\Chev$ over the ring of Laurent
    polynomials. He shows that this group has an {\small RGD}-system
    (confer~\cite[Definition~2, pages~14f]{Abramenko:1996}).

    Finally one can derive the Moufang axioms from the
    {\small RGD}-axioms. This is not too difficult since the Moufang
    axioms~\ref{MoufangZwei} to~\ref{MoufangVier} can be read
    as geometric interpretations of analogous {\small RGD}-axioms.
    The transitivity of the action in the Moufang axiom~\ref{MoufangEins}
    follows from~\cite[5.6~Proposition~3, page~564]{Tits:1987} whereas
    ({\small RGD}~3) immediately implies that the action is
    simply-transitive.
  \end{sketch}

  In~\cite{Aberg:1986}, H \AA{}berg gave a method to detect the
  vanishing of homotopy groups for certain subspaces in products of
  trees. We will generalize his ideas to affine buildings.
  \begin{Def}
    Let $\Build$ be an affine building and
    $\InfChamber$ a chamber at infinity. We call a
    (not necessarily infinite) sequence of apartments
    $\FamOf{\ApartmentInd{\Index}}{}$ 
    \notion{directed} if for each index~$\SupIndex$ the closed set
    $$
      \Closure{%
        \ApartmentInd{\SupIndex}
        \setminus
        \Union{\SubIndex<\SupIndex}{\ApartmentInd{\SubIndex}}
      }
    $$
    is an intersection of half apartments in $\ApartmentInd{\SupIndex}$
    that do not contain $\InfChamber$.
  \end{Def}
  \begin{prop}\label{PositiveRichtung.Haupthilfslemma}
    If $\Build$ is locally finite and Moufang then there is an infinite
    directed sequence
    $\FamOf{\ApartmentInd{\Index}}{\Index\in\NNN}$
    of apartments that covers $\Build$.
  \end{prop}
  The proof will take the remainder of this section.
  So let us fix $\Build$ and $\InfChamber$ as in the proposition. Further let us
  consider a bi-infinite geodesic gallery
  \begin{notationlist}
    \item
      $\gallery:=\Gallery{\ldots,\AffChamberInd{\Index},\AffChamberInd{\Index+1},
      \ldots}$
      within a chosen standard apartment
    \item
      $\AffApp$ that contains $\InfChamber$. We will make more specific choices later. Let
    \item
      $\AltAffRootInd{\Index}$ be the half apartment of
      $\AffApp$ containing $\AffChamberInd{\Index-1}$ but not
      $\AffChamberInd{\Index}$. We will write
    \item
      $\RootGroupTo{\Index}:=\RootGroupTo{\AltAffRootInd{\Index}}$
      to avoid double subscripts.
  \end{notationlist}
  The following facts are quoted from~\cite[pages~74+]{Ronan:1989}
  or can be proved by similar arguments without difficulty:
  \begin{enumerate}
    \item\label{Moufang:FixiertStern}
      $\RootGroupTo{\AffRoot}$ fixes the star of every panel
      in $\AffRoot \setminus \Boundary{\AffRoot}$.
    \item\label{Moufang:RootGroupProd}
      $
        \RootGroupProdFromTo{\LowInd}{\HighInd}
        :=
        \GroupGen{
          \RootGroupTo{\Index}
          \suchthatvrule
          \Index \in \SetOf{\LowInd,\ldots,\HighInd}
        }
        =
        \RootGroupTo{\LowInd} \cdots \RootGroupTo{\HighInd}
      $.
      Moreover, the corresponding factorization
      $
        \RootGroupElement
        =
        \RootGroupElementTo{\LowInd} \cdots \RootGroupElementTo{\HighInd}
      $
      is unique for every element
      $
        \RootGroupElement
        \in
        \RootGroupProdFromTo{\LowInd}{\HighInd}
      $.
    \item\label{Moufang:Normalisiert}
      $\RootGroupTo{\LowInd-1}$ and
      $\RootGroupTo{\HighInd+1}$ normalize
      $\RootGroupProdFromTo{\LowInd}{\HighInd}$.
    \item\label{Moufang:RootGroupProdEndlich}
      $\RootGroupProdFromTo{\LowInd}{\HighInd}$ is finite
      since the building $\Build$ is locally finite.
  \end{enumerate}
  \begin{Lemma}\label{Apartment.Durchschnitt}
    An element
    $
      \RootGroupElement
      =
      \RootGroupElementTo{\LowInd} \cdots \RootGroupElementTo{\HighInd}
      \in
      \RootGroupProdFromTo{\LowInd}{\HighInd}
    $
    fixes exactly those chambers of $\AffApp$ that lie in the intersection
    of all $\AltAffRootInd{\Index}$ with $\RootGroupElementTo{\Index}\neq1$.
    Moreover, this is the intersection
    $
      \RootGroupElement \AffApp
      \cap
      \AffApp
    $.
  \end{Lemma}
  \begin{proof}
    Put
    $
      \RootGroupElementUpTo{\Index}
      :=
      \RootGroupElementTo{\LowInd} \cdots \RootGroupElementTo{\Index}
    $.
    We proceed by induction on $\Index$.
    For $\Index = \LowInd$, we have
    $\RootGroupElementUpTo{\Index}\in\RootGroupTo{\LowInd}$.
    If this element is trivial, then it fixes all of $\AffApp$, which we regard
    as the intersection of an empty family of half apartments.
    Now suppose $\RootGroupElementTo{\LowInd}\neq 1$ and assume that there was
    a chamber $\AffChamber$
    in $\AffApp\setminus\AltAffRootInd{\LowInd}$ fixed by
    $\RootGroupElementUpTo{\LowInd}$. Then the element
    $\RootGroupElementUpTo{\LowInd}$ would fix every minimal gallery connecting
    $\AffChamber$ to a chamber in $\AltAffRootInd{\LowInd}$. In this case
    there would be a panel $\Panel\in\Boundary{\AltAffRootInd{\LowInd}}$
    whose star contained two chambers fixed by
    $\RootGroupElementUpTo{\LowInd}$. This, however, is impossible
    by axiom~\ref{MoufangEins}.

    Now suppose $\Index > \LowInd$. By induction, the lemma holds for
    $\RootGroupElementUpTo{\Index-1}$. Put
    $
      \IntersectionOne
      :=
      \RootGroupElementUpTo{\Index-1} \AffApp
      \cap
      \AffApp
    $.
    We have to show that for $\RootGroupElementTo{\Index} \neq 1$,
    the identity
    $$
      \LeftPart
      :=
      \IntersectionOne \cap \AltAffRootInd{\Index}
      =
      \RootGroupElementUpTo{\Index} \AffApp \cap \AffApp
      =:
      \RightPart
    $$
    holds. By induction, we see that
    $\AffChamberInd{\LowInd-1} \in \AltAffRootInd{\LowInd}
     \cap \cdots \cap \AltAffRootInd{\Index}
     \subseteq \LeftPart \subseteq \RightPart
    $
    since $\RootGroupElementTo{\Index}$ fixes the half apartment
    $\AltAffRootInd{\Index}$.

    Assume there is a chamber
    $\AffChamber \in \RightPart \setminus \LeftPart$.
    Observe that
    $\RootGroupElementUpTo{\Index}$ stabilizes $\AffChamberInd{\LowInd-1}$
    and preserves the type of minimal galleries. The chamber
    $\AffChamber$ is uniquely determined by the type of a minimal
    gallery in $\AffApp$ connecting it to
    $\AffChamberInd{\LowInd-1}$.
    Hence $\RootGroupElementUpTo{\Index}$ fixes $\AffChamber$
    and also the gallery connecting $\AffChamber$
    to $\AffChamberInd{\LowInd-1}$.
    Consider the panel $\Panel\in\Boundary{\LeftPart}$ within this gallery.
    If $\Panel \in \Boundary{\AltAffRootInd{\Index}}$ then both of its
    neighboring chambers in $\AffApp$ lie in $\IntersectionOne$ and
    $\RootGroupElementTo{\Index}$ fixes both of them which
    cannot happen. Therefore, $\Panel$ lies in the interior
    of $\AltAffRootInd{\Index}$. In this case,
    $\RootGroupElementTo{\Index}$ fixes the star of the panel $\Panel$
    whence $\RootGroupElementUpTo{\Index-1}$ has to fix both neighboring
    chambers in $\AffApp$. However, one of them does
    not belong to $\IntersectionOne$ whence
    $\RootGroupElementUpTo{\Index-1}$ does not fix it.
    Thus, the hypothesis $\LeftPart\subsetneqq\RightPart$
    leads to a contradiction.
  \end{proof}

  We will construct a sequence of automorphisms
  $\AutomInd{\SupIndex}\in
  \GroupGen{\RootGroupTo{\Index}\suchthatvrule\Index\in\ZZZ}
  \subseteq\AutOf{\Build}$ that will give rise to a
  directed sequence
  $\ApartmentInd{\SupIndex}:=\AutomIndEv{\SupIndex}{\AffApp}$.
  So we call a possibly finite sequence
  $\AutomInd{\One},\ldots,\AutomInd{\Index},\ldots$ of automorphisms
  of $\Build$ \notion{directed} if for each index $\SupIndex$, the set
  $$
    \Union{\SubIndex<\SupIndex}{%
      \Ev{\AutomInd{\SupIndex}^{-1}\AutomInd{\SubIndex}}{\AffApp}
      \cap \AffApp
    }
  $$
  is \notion{coconvex}, ie, the complement of a convex set.  Here, we
  use the notion of convexity that stems from the Euclidean metric on
  $\AffApp$. However, recall that a subcomplex of the Coxeter complex
  $\AffApp$ is (metrically) convex if and only if it is an
  intersection of half apartments.
  \begin{observation}\label{directed} Suppose the gallery~$\gallery$
  is chosen such that every half apartment $\AltAffRootInd{\Index}$
  contains the chamber $\InfChamber$.  Then every intersection \(
  \Ev{\AutomInd{\SupIndex}^{-1}\AutomInd{\SubIndex}}{\AffApp} \cap
  \AffApp \) contains $\InfChamber$.  In this case, for any directed
  sequence $\FamOf{\AutomInd{\Index}}{\Index\in\NNN}$ of automorphisms
  of $\Build$, the sequence
  $\FamOf{\AutomIndEv{\Index}{\AffApp}}{\Index\in\NNN}$ is a directed
  sequence of apartments.\qed \end{observation}
  \begin{Lemma}\label{Automorphismen.Aufzaehlung} Every directed
  enumeration of $\RootGroupProdFromTo{\LowInd}{\HighInd}$ can be
  extended to a directed enumeration of
  $\RootGroupProdFromTo{\LowInd-1}{\HighInd+1}$.  \end{Lemma}
  \begin{proof} We will only show that we can extend to
  $\RootGroupProdFromTo{\LowInd}{\HighInd+1}$. The argument for the
  lower index is similar.

    Let
    $
     \RootGroupProdFromTo{\LowInd}{\HighInd}=
     \SetOf{\AutomInd{\One},\ldots,\AutomInd{\BlockTop}}
    $
    be a directed sequence enumerating
    $\RootGroupProdFromTo{\LowInd}{\HighInd}$ and choose
    an enumeration
    $
     \RootGroupAt{\HighInd+1}=
     \SetOf{\RootGroupElementTo{\One},\ldots,\RootGroupElementTo{\GroupTop}}
    $
    starting with $\RootGroupElementTo{\One}=1$.
    Then
    $$
      \SetOf{%
\RootGroupElementTo{\One}\AutomInd{\One},\ldots,\RootGroupElementTo{\One}\AutomInd{\BlockTop}
,\ldots,
\RootGroupElementTo{\Two}\AutomInd{\One},\ldots,\RootGroupElementTo{\Two}\AutomInd{\BlockTop}
,\ldots,
\RootGroupElementTo{\GroupTop}\AutomInd{\One},\ldots,\RootGroupElementTo{\GroupTop}\AutomInd{\BlockTop}
      }
    $$
    is an enumeration of
    $\RootGroupProdFromTo{\LowInd}{\HighInd+1}$ which splits into $\GroupTop$
    blocks of length $\BlockTop$. The first block is the given enumeration
    of $\RootGroupProdFromTo{\LowInd}{\HighInd}$. So we extended
    the enumeration. We will prove that this extension is directed.
    Within a block, quotients of elements agree with quotients in the
    first block. Thus, it suffices to prove that
    sets of the form
    $$
      \biggl(
        \Union{%
          \substack{%
            \GroupSubInd<\GroupSupInd\\
            \BlockSubInd\leq\BlockTop
          }
        }{%
          \Ev{%
            \AutomInd{\BlockSupInd}^{-1}
            \RootGroupElementTo{\GroupSupInd}^{-1}\RootGroupElementTo{\GroupSubInd}
            \AutomInd{\BlockSubInd}
          }{%
            \AffApp
          }
          \cap
          \AffApp
        }
      \biggr)
      \cup
      \biggl(
        \Union{%
          \BlockSubInd<\BlockSupInd
        }{%
          \Ev{%
            \AutomInd{\BlockSupInd}^{-1}\AutomInd{\BlockSubInd}
          }{%
            \AffApp
          }
          \cap
          \AffApp
        }
      \biggr)
    $$
    are coconvex.

    The right hand term is coconvex because we extended
    a directed sequence.

    If the left hand side is not empty it equals
    the half apartment $\AltAffRootInd{\HighInd+1}$:
    We have the factorization
    $$
      \AutomInd{\BlockSupInd}^{-1}
      \RootGroupElementTo{\GroupSupInd}^{-1}\RootGroupElementTo{\GroupSubInd}
      \AutomInd{\BlockSubInd}
      =
      \underbrace{
      (
        \AutomInd{\BlockSupInd}^{-1}
        \RootGroupElementTo{\GroupSupInd}^{-1}\RootGroupElementTo{\GroupSubInd}
        \AutomInd{\BlockSubInd}
        \RootGroupElementTo{\GroupSubInd}^{-1}\RootGroupElementTo{\GroupSupInd}
      )}_{
        \in\RootGroupProdFromTo{\LowInd}{\HighInd}
      }
      \underbrace{(
        \RootGroupElementTo{\GroupSupInd}^{-1}\RootGroupElementTo{\GroupSubInd}
      )}_{
        \in\RootGroupAt{\HighInd+1}
      } .
    $$
    Since $\RootGroupElementTo{\GroupSupInd}^{-1}\RootGroupElementTo{\GroupSubInd}\neq1$,
    Lemma~\ref{Apartment.Durchschnitt} implies
    $
          \Ev{%
            \AutomInd{\BlockSupInd}^{-1}
            \RootGroupElementTo{\GroupSupInd}^{-1}\RootGroupElementTo{\GroupSubInd}
            \AutomInd{\BlockSubInd}
          }{%
            \AffApp
          }
          \cap
          \AffApp
      \subseteq
      \AltAffRootInd{\HighInd+1}
      .
    $
    As $\BlockSubInd\in\SetOf{\One,\ldots,\BlockTop}$ varies, 
    the first factor
    $
      \AutomInd{\BlockSupInd}^{-1}
      \RootGroupElementTo{\GroupSupInd}^{-1}\RootGroupElementTo{\GroupSubInd}
      \AutomInd{\BlockSubInd}
      \RootGroupElementTo{\GroupSubInd}^{-1}\RootGroupElementTo{\GroupSupInd}
    $
    runs through all of $\RootGroupProdFromTo{\LowInd}{\HighInd}$. In particular,
    at some point, $\AutomInd{\BlockSupInd}^{-1}
      \RootGroupElementTo{\GroupSupInd}^{-1}\RootGroupElementTo{\GroupSubInd}
      \AutomInd{\BlockSubInd}
      \RootGroupElementTo{\GroupSubInd}^{-1}\RootGroupElementTo{\GroupSupInd}$ will be the trivial element.
    Therefore,
    $\AltAffRootInd{\HighInd+1}$ contains all sets
    $
        \Ev{%
          \AutomInd{\BlockSupInd}^{-1}
          \RootGroupElementTo{\GroupSupInd}^{-1}\RootGroupElementTo{\GroupSubInd}
          \AutomInd{\BlockSubInd}
        }{%
          \AffApp
        }
        \cap
        \AffApp
    $
    and equals at least one of them.
  \end{proof}

  Now, we choose the gallery~$\gallery$.
  We start
  with a straight line in the affine space $\AffApp$ that
  intersects each sector representing $\InfChamber$. Note that
  any such line also intersects each sector representing the opposite
  chamber at infinity: geometrically such a line connects two
  antipodal points (in the spherical metric) of the spherical building
  at infinity.
  We move this line in its parallelity class so that it does
  not intersect the codimension~$2$~skeleton of $\AffApp$. We can
  do this, as the set of forbidden lines has measure $0$ since the
  space of parallel lines has codimension~$1$.
  Given such a line, it intersects a sequence of
  chambers and thus defines the
  bi-infinite gallery
  \(
    \gallery
    =
    \Gallery{
      \ldots,
      \AffChamberInd{\Index},
      \AffChamberInd{\Index+1},
      \ldots
    }
    .
  \)
  By construction, all
  half apartments~$\AltAffRootInd{\Index}$ contain the
  chamber~$\InfChamber$, and the apartment~$\AffApp$ is
  the convex hull of $\gallery$, ie, there is no half apartment
  containing the whole gallery $\gallery$.
  
  By Lemma~\ref{Automorphismen.Aufzaehlung}, we can find a directed
  sequence enumerating the group
  $\RootGroupProdFromTo{-\infty}{\infty}
   :=\GroupGen{\RootGroupTo{\Index}\suchthatvrule\Index\in\ZZZ}
   =\Union{\LowInd,\HighInd,\LowInd \leq \HighInd}{
     \RootGroupTo{\LowInd} \cdots \RootGroupTo{\HighInd}
    }
   =\SetOf{\AutomInd{\Index}\suchthatvrule\Index\in\NNN}
  $.
  By Observation~\ref{directed},
  this sequence induces a directed sequence
  $\ApartmentInd{\Index}:=\AutomIndEv{\Index}{\AffApp}$
  of apartments.
  To finish the proof of Proposition~\ref{PositiveRichtung.Haupthilfslemma},
  it remains to show that these apartments cover $\Build$.

  Let $\AltChamber$ be a chamber in $\Build$. It suffices to show
  that there is an element $\Autom\in\RootGroupProdFromTo{-\infty}{\infty}$
  such that $\Autom\AltChamber\subset\AffApp$. We will use induction
  on the length of a shortest gallery
  \(
      \AltChamber=\AltChamberInd{0},
      \AltChamberInd{1},
      \ldots,
      \AltChamberInd{\TopIndex-1},
      \AltChamberInd{\TopIndex}\subset\AffApp
  \).
  Let $\root$ be the half apartment that contains $\InfChamber$ and
  whose boundary contains the panel
  $\AltChamberInd{\TopIndex-1}\cap\AltChamberInd{\TopIndex}$.
  Since no half apartment contains the whole gallery $\gallery$, we
  have
  \(
    \RootGroupAt{\root}
    \subgroup
    \RootGroupProdFromTo{-\infty}{\infty}
    .
  \)
  By Axiom~(\ref{MoufangEins}), there is an element
  $\AutomInd{\root}\in\RootGroupAt{\root}$ such that
  $\AutomInd{\root}\AltChamberInd{\TopIndex-1}\subset\AffApp$.
  Hence, $\AutomInd{\root}\AltChamber$ can be connected
  to $\AffApp$ by a shorter gallery than $\AltChamber$.
  By induction, there is an element
  $\Autom'\in\RootGroupProdFromTo{-\infty}{\infty}$ such that
  $\Autom'\AutomInd{\root}\AltChamber\subset\AffApp$.
  This completes the proof of
  Proposition~\ref{PositiveRichtung.Haupthilfslemma}.\qed

\section{Tameness, connectivity, and finiteness properties}
\label{EndProof}
  Let
      $\SheetCompl$ be a metric combinatorial \CW-complex such
      that the restriction of the metric to a closed
      cell makes the cell isometric to a convex Euclidean polyhedron, and let
      $\SheetProj \mapcolon  \SheetCompl \rightarrow \Eukl$ be a projection
      onto a Euclidean space
      $\Eukl$. Furthermore, let us fix a non-empty set
      $\SheetSet$ of linear forms on $\Eukl$. 
  A \notion{sheet} is a subcomplex of $\SheetCompl$ that $\SheetProj$
  maps isometrically onto $\Eukl$. A closed subset of $\Eukl$ is
  \notion{$\SheetSet$--convex} if it is an
  intersection of halfspaces
  \[
    \Intersection{\SheetLin\in\SheetSet'}{
      \SetOf{\EuklPoint\in\Eukl \suchthatvrule \SheetLinEv{\EuklPoint}\geq\SheetConstAt{\SheetLin}}
    }
  \]
  where $\SheetSet'$ is a non-empty subset of $\SheetSet$ and $\SheetConstAt{\SheetLin}\in\RRR$ for
  every $\SheetLin\in\SheetSet'$. In particular, $\Eukl$ is \emph{not} $\SheetSet$--convex unless
  $0\in\SheetSet$. A closed subset of a sheet is called
  $\SheetSet$--convex if it is identified with a $\SheetSet$--convex subset
  of $\Eukl$ via $\SheetProj$. We follow the terminology
  of~\cite{Bux:1997} and call the triple
  $\TupelOf{\SheetCompl, \Eukl, \SheetProj}$ a
  \notion{$\SheetSet$--complex} if every compact subset of $\SheetCompl$
  is covered by an \notion{increasing sequence} of sheets.
  Here we call a sequence $\FamOf{\SheetInd{\Index}}{\Index}$
  increasing if, at every stage, the piece
  $\Closure{
     \SheetInd{\HighInd}
     \setminus
     \Union{\LowInd<\HighInd}{\SheetInd{\LowInd}}
   }$
  just added is $\SheetSet$--convex.
  Note that the set of sheets in a $\SheetSet$--complex
  $\SheetCompl$ form a cover of $\SheetCompl$.

  \begin{observation}\label{GebaeudeAtPlace.SheetComplex}
    Proposition~\ref{PositiveRichtung.Haupthilfslemma} applies in particular
    to the building $\AffBuildAt{\place}$ and the chamber
    $\InfChamberAt{\place}$. Thus,
    $\TupelOf{\AffBuildAt{\place}, \AffAppAt{\place}, \ProAt{\place}}$
    is a $\NegativeRootsAt{\place}$--complex.\qed
  \end{observation}
  By \cite[Example~7.2]{Bux:1997}, we can lift this result to the product $\AffAppProd$:
  \begin{cons}
    Each linear form in $\NegativeRootsAt{\place}$ induces
    by composition with the projection
    $\AffAppProd\rightarrow\AffAppAt{\place}$
    a linear form on $\AffAppProd$. In this sense, put
    \begin{notationlist}
      \item
        $\NegativeRootsUnion:=
        \DisjointUnion{\place\in\Set}{\NegativeRootsAt{\place}}$.
    \end{notationlist}
    Then
    $\TupelOf{\AffBuildProd,\AffAppProd,\ProProd \mapcolon \AffBuildProd\rightarrow\AffAppProd}$
    is a $\NegativeRootsUnion$--complex.\qed
  \end{cons}
   
  We will study the action of $\BorelOver{\OkaSet}$ on the complex
  \begin{notationlist}
    \item
      $\BorelSpace := \PreImOf{\ProSum}{0} = \PreImOf{\ProProd}{\Kernel}$
      where
    \item
      $\Kernel := \Kern{\CoordinatesSum}$. The map
      $\ProProd \mapcolon \AffBuildProd \rightarrow \AffAppProd$ obviously
      restricts to $\BorelSpace$ and yields a projection
      $\BorelSpace \rightarrow \Kernel$, which we will also denote by
      $\ProProd$.

      The linear forms in $\NegativeRootsUnion$ restrict to $\Kernel$. Let
    \item
      $\NegativeRootsRestricted := \Restr{\NegativeRootsUnion}{\Kernel}$
      be the system of linear forms obtained this way.
  \end{notationlist}
  \begin{observation}\label{BorelSpace}
    Let $\TupelOf{\SheetCompl,\Eukl,\SheetProj}$ be a $\SheetSet$--complex
    and $\EuKern$ be a subspace of $\Eukl$. Then
    $\TupelOf{
      \PreImOf{\SheetProj}{\EuKern},\EuKern,\Restr{\SheetProj}{\EuKern}
    }$ is a $\Restr{\SheetSet}{\EuKern}$--complex.
    Here, $\Restr{\SheetSet}{\EuKern}$ denotes the set of
    those linear forms in $\EuKern$ obtained by restricting a
    linear form in $\SheetSet$ to $\EuKern$.
    In particular, $\TupelOf{\BorelSpace,\Kernel,\ProProd}$ is a
    $\NegativeRootsRestricted$--complex.\qed
  \end{observation}
    Note that the metric on $\PreImOf{\SheetProj}{\EuKern}$
    is not the path metric but the restriction of the metric on
    $\SheetCompl$. Moreover, in general, elements of $\SheetSet$ could
    vanish on $\EuKern$.
    Observe, however, that $0\not\in\NegativeRootsRestricted$ because of
    the product formula.
  \begin{lemma}\label{Urbild:Kernel:Zusammenhaengend}
    $\BorelSpace$ is $(\CardOf{\Set}-2)$--connected.
  \end{lemma}
  \begin{proof}
    By~\cite[Lemma~7.3]{Bux:1997} 
    $\BorelSpace=\ProSum^{-1}(0)=\ProProd^{-1}(\Kernel)$ is
    $(\CardOf{\Set}-2)$--connected provided
    $\NegativeRootsRestricted$ is \notion{$(\CardOf{\Set}-1)$--tame},
    ie, no positive linear combination of up to
    $\CardOf{\Set}-1$ elements in $\NegativeRootsRestricted$
    vanishes.

    Let us call a positive combination of forms in
    $\NegativeRootsUnion=
    \DisjointUnion{\place\in\Set}{\NegativeRootsAt{\place}}$
    \xnotion{complete} if it involves for each $\place\in\Set$ at
    least one element of $\NegativeRootsAt{\place}$ non-trivially.
    We shall argue that any positive combination of forms in
    $\NegativeRootsUnion$ that vanishes on $\Kernel$ is complete
    and hence involves at least $\CardOf{\Set}$ forms with strictly
    positive weight. This implies that
    $\NegativeRootsRestricted=\Restr{\NegativeRootsUnion}{\Kernel}$
    is $(\CardOf{\Set}-1)$--tame.

    Note that a positive combination of complete forms is complete. Since
    each element of $\NegativeRootsAt{\place}$ is a positive combination
    of of base roots
    $\bTupelOf{\AffRootAtInd{\place}{1},\ldots,\AffRootAtInd{\place}{\TorDim}}$,
    we can write any positive combination of elements of
    $\NegativeRootsUnion$ as a positive combination of base roots.
    Thus, it suffices to prove:
    \begin{quote}
    If
    $
      \Sum{\Index,\place}{
        \KoeffizientAtInd{\place}{\Index}\AffRootAtInd{\place}{\Index}
      }
    $
    is a positive combination of base roots that vanishes on $\Kernel$, then
    for each $\place\in\Set$, there is an $\Index\in\SetOf{1,\ldots,\TorDim}$
    with $\KoeffizientAtInd{\place}{\Index}>0$.
    \end{quote}
    This, however, follows from the definition of $\Kernel$ since
    $
      \Sum{\Index,\place}{
        \KoeffizientAtInd{\place}{\Index}\AffRootAtInd{\place}{\Index}
      }
    $
    vanishes on $\Kernel$ only if, for each $\Index$, the coefficient
    $\KoeffizientAtInd{\place}{\Index}$ does not depend on $\place$.
  \end{proof}
  We are now in a position to tackle the companion of
  Theorem~\ref{ObereAbschaetzung}.
  \begin{theorem}\label{UntereAbschaetzung}
    The group $\BorelOver{\OkaSet}$ is of type \F{\CardOf{\Set}-1}.
  \end{theorem}
  \begin{proof}
    $\BorelOver{\OkaSet}$ acts on $\BorelSpace$ by cell-permuting homeomorphisms.
    The action is cocompact by
    Lemma~\ref{Operation:kokompakt}. The cell stabilizers are finite
    by Lemma~\ref{Stabilisatoren:endlich}. The space $\BorelSpace$ is
    $(\CardOf{\Set}-2)$--connected by Lemma~\ref{Urbild:Kernel:Zusammenhaengend}.
    Hence, the claim follows from Corollary~\ref{Brown:Simple_Criterion}.
  \end{proof}
  \begin{rem}
    The proof of Theorem~\ref{UntereAbschaetzung} mimics
    the argument given in~\cite{Bux:1997} and uses the fact that
    the buildings $\AffBuildAt{\place}$ are Moufang, which is
    a global property.
    Hence, this method does not apply to arithmetic groups over
    number fields as the associated buildings are not Moufang.
    However, the discussion in Section~\ref{Sec:SlZwei}
    indicates that it might be possible to find an argument
    for Theorem~\ref{UntereAbschaetzung}
    somewhat in the spirit~\cite{Bestvina.Brady:1997} that
    relies only on local properties of the buildings
    $\AffBuildAt{\place}$.
  \end{rem}
  Theorems~\ref{ObereAbschaetzung} and~\ref{UntereAbschaetzung} together
  imply Theorem~A stated in the introduction.\qed

\section{The geometric invariants of Bieri, Geoghegan, Neumann, Strebel and Renz}
  \label{Sec:Sigma}
  ``Geometric invariants'', born in~\cite{Bieri.Strebel:1980}, are by now
  connected with five names. It would be unfeasible to draw
  a complete picture of the theory or its history.
  The reader may consult~\cite{Renz:1988} or \cite{Bieri.Renz:1988}
  where high dimensional invariants were first introduced.
  A more recent source is~\cite{Bieri:1999}.

  The invariant $\BNS{\type}{\Group}$ is defined for every
  group~$\Group$ of type \F{\type}. It is an open cone
  in $\GrHom{\Group}{\RRR}$. For each group $\Group$, we have a chain
  of inclusions
  \[
    \GrHom{\Group}{\RRR} \supseteq
    \BNS{1}{\Group} \supseteq \BNS{2}{\Group} \supseteq \dots \supseteq
    \BNS{\type}{\Group}
    .
  \]
  
  \begin{Def}
    Let $\Group$ be a group of type \F{\type} and $\Komplex$ be an
    $(\type-1)$--connected \CW-complex upon which $\Group$ acts by
    cell permuting homeomorphisms such that the stabilizer of
    each $\celldim$--cell is of type \F{\type-\celldim} and such that
    the action on the $\type$--skeleton of $\Komplex$ is cocompact.
    For each non-trivial homomorphism $\homo \in \GrHom{\Group}{\RRR}$,
    there exists a $\Group$--equivariant \notion{height function}
    on $\Komplex$, ie, a continuous map
    $\Height{\homo} \mapcolon \Komplex \rightarrow \RRR$ satisfying
    $\HeightAt{\homo}{\GroupElement\KomplexPoint} =
    \homoat{\GroupElement} + \HeightAt{\homo}{\KomplexPoint}$
    for all points $\KomplexPoint\in\Komplex$ and all elements
    $\GroupElement\in\Group$.

    Let $\Height{\homo}$ be a height function associated to $\homo$.
    By definition $\homo \in \BNS{\type}{\Group}$ holds if and only if
    $\Komplex$ is \notion{essentially $(\type-1)$--connected
    with respect to $\Height{\homo}$}, ie, the directed systems
    $\FamOf{%
       \FundamentalGroup{\Index}{%
         \PreImOf{\Height{\homo}}{[\real,\infty)}
       }
     }%
     {\real\in\RRR}$
    of homotopy groups are essentially trivial
    for all $\Index < \type$. For the definition of an essentially trivial
    directed system of groups, see Citation~\ref{Kriterium:Brown}.
  \end{Def}
  \begin{rem}
    One can show that this definition is independent of the complex $\Komplex$ and the
    height function $\Height{\homo}$ chosen to represent $\homo$.
    See  \cite[Theorem~12.1]{Bieri.Geoghegan:2002}.
    In earlier stages of the theory, the action of $\Group$ on $\Komplex$ was
    required to be free. The weaker hypothesis regarding finiteness properties of
    cell stabilizers was proven to yield an equivalent definition of $\BNS{\type}{\Group}$
    by H Meinert in his PhD thesis. Its content is published in
    his papers~\cite{Meinert:1996} and~\cite{Meinert:1997}.
    However, the reader might prefer the reference~\cite[Theorem~12.1]{Bieri.Geoghegan:2002}
    given above since that approach
    avoids a detour via homological invariants.
  \end{rem}
  \begin{rem}\label{Sufficient}
    Still, our definition differs slightly from the literature. Since the difference
    lies in the order of quantifiers, this might be one of the few instances where
    logical formalism actually eases understanding.
    In the definition above, we require:
    \[
      \mathop{\forall}\limits_{\Index<\type}\,
      \mathop{\forall}\limits_{\HighReal\in\RRR}\,
      \mathop{\exists}\limits_{\LowReal\leq\HighReal}\,
      \FundamentalGroup{\Index}{%
        \PreImOf{\Height{\homo}}{[\HighReal,\infty)}
      }
      \rightarrow
      \FundamentalGroup{\Index}{%
        \PreImOf{\Height{\homo}}{[\LowReal,\infty)}
        }
      \text{\ vanishes.}
    \]
    Usually, the condition is phrased in a slightly stronger way:
    \[
      \mathop{\forall}\limits_{\Index<\type}\,
      \mathop{\exists}\limits_{\TheLag\geq0}\,
      \mathop{\forall}\limits_{\HighReal\in\RRR}\,
      \FundamentalGroup{\Index}{%
        \PreImOf{\Height{\homo}}{[\HighReal,\infty)}
      }
      \rightarrow
      \FundamentalGroup{\Index}{%
        \PreImOf{\Height{\homo}}{[\HighReal-\TheLag,\infty)}
        }
      \text{\ vanishes.}
    \]
    However, since $\Group$ acts non-trivially on $\RRR$ by translations,
    both definitions are equivalent
    to the following even weaker condition:
    \[
      \mathop{\forall}\limits_{\Index<\type}\,
      \mathop{\exists}\limits_{\HighReal\in\RRR}\,
      \mathop{\exists}\limits_{\LowReal\leq\HighReal}\,
      \FundamentalGroup{\Index}{%
        \PreImOf{\Height{\homo}}{[\HighReal,\infty)}
      }
      \rightarrow
      \FundamentalGroup{\Index}{%
        \PreImOf{\Height{\homo}}{[\LowReal,\infty)}
        }
      \text{\ vanishes.}
    \]
    In particular, if $\PreImOf{\Height{\homo}}{[\HighReal,\infty)}$ is $(\type-1)$--connected
    for some $\HighReal$, we have $\character\in\BNS{\type}{\Group}$.
  \end{rem}
  Generally, the geometric invariants of a group $\Group$ are
  hard to compute---they convey a lot of information. Eg,
  $\BNS{\type}{\Group}$ determines for each normal subgroup
  $\SubGroup\normalsubgroup\Group$ with abelian factor whether $\SubGroup$ is of type \F{\type}:
  \begin{cit}[{confer \cite[Satz~C, page~17]{Renz:1988}}]%
    \label{Kriterium:Renz}
    Let $\Group$ be a group of type \F{\type} and
    $\SubGroup\normalsubgroup\Group$ be a normal subgroup with abelian factor group.
    Then, $\SubGroup$ is of type \F{\type} if and only if
    all non-trivial homomorphisms in $\GrHom{\Group}{\RRR}$
    that vanish on $\SubGroup$ belong to $\BNS{\type}{\Group}$.
  \end{cit}

  We already know that $\BorelOver{\OkaSet}$ is of type \F{\CardOf{\Set}-1}.
  Thus, the geometric invariants $\BNS{\type}{\BorelOver{\OkaSet}}$
  are defined for $1\leq\type<\CardOf{\Set}$. Since we are working in
  positive characteristic, all unipotent groups are torsion groups.
  Thus, every homomorphism from $\BorelOver{\OkaSet}$ to $\RRR$ factors via
  $\TorusOver{\OkaSet}$ which spans a maximal lattice in the vector
  space $\Kernel$ that can, therefore, be described as
  $\Kernel \isom \Tensor{\ZZZ}{\TorusOver{\OkaSet}}{\RRR}$.
  Thus,
  \begin{notationlist}
    \item
      $\KernDual :=
       \Homs{\RRR}{\Kernel}{\RRR}\isom
       \Homs{\ZZZ}{\TorusOver{\OkaSet}}{\RRR}\isom
       \GrHom{\BorelOver{\OkaSet}}{\RRR}\supseteq
       \BNS{\type}{\BorelOver{\OkaSet}}
      $.
      To avoid cumbersome notation, we will frequently use the complement
    \item
      $\BNSc{\type}{\BorelOver{\OkaSet}}
       := \KernDual \setminus \BNS{\type}{\BorelOver{\OkaSet}}
      $.
  \end{notationlist}

  We will use the action of $\BorelOver{\OkaSet}$ on the space
  $\BorelSpace$ defined in Section~\ref{EndProof} to
  compute the geometric invariants of $\BorelOver{\OkaSet}$.
  A homomorphism from $\TorusOver{\OkaSet}$ to the real numbers $\RRR$ can be
  represented by a linear form on $\Kernel$. Composition with
  $\ProProd \mapcolon  \BorelSpace \rightarrow \Kernel$ yields a
  $\BorelOver{\OkaSet}$--equivariant height function on $\BorelSpace$.

  We will derive a lower and an upper bound for
  $\BNSc{\type}{\BorelOver{\OkaSet}}$. These bounds will be sharp
  only in the rank--$1$--case. To state these bounds succinctly, we
  need some notation. Let
  \begin{notationlist}
    \item
      $\BaseRootsUnion :=
       \subseteq
       \NegativeRootsUnion
      $
      denote the set of all linear forms on $\AffAppProd$ induced by
      base roots and let
    \item
      $\BaseRootsRestricted := \Restr{\BaseRootsUnion}{\Kernel}
      \subseteq\NegativeRootsRestricted$ 
      be the restriction of $\BaseRootsUnion$ to $\Kernel$. By
    \item
      $\ConvSub{\type}{\Set}{\BaseRootsRestricted}$ we refer to
      the set of all linear forms on $\Kernel$ that can be
      written as a positive linear combination
      $\Sum{\place\in\Set}{%
         \KoeffizientAt{\place}
         \AffRootAtInd{\place}{\IndexAt{\place}}
       }$
      wherein at most $\type$ coefficients do not vanish. Note that
      all base roots occurring in the combination belong to distinct
      places.

      For an arbitrary set $\SigmaSet$ of linear forms, let
    \item
      $\Conv{\type}{\SigmaSet}$ be the set of all linear forms
      that are positive combinations of up to $\type$ elements
      of $\SigmaSet$.
  \end{notationlist}
  \begin{theorem}\label{Sigma}
    For $\type < \CardOf{\Set}$, we have the chain of inclusions
    \[
      \ConvSub{\type}{\Set}{\BaseRootsRestricted}
      \subseteq\BNSc{\type}{\BorelOver{\OkaSet}}
      \subseteq
      \Conv{\type}{\NegativeRootsRestricted}
      .
    \]
  \end{theorem}
  \begin{rem}
    In the rank--$1$--case, the upper and the lower bound coincide.
    Hence they determine $\BNSc{\type}{\BorelOver{\OkaSet}}$.
    In this case, $\BorelOver{\OkaSet}$ is metabelian, and
    Theorem~\ref{Sigma} establishes the $\Bns{\type}$--conjecture for
    $\BorelOver{\OkaSet}$.
    This also follows from results of
    D Kochloukova \cite{Kochloukova:1997}.
  \end{rem}
  \begin{rem}
    It is easy to spell out some consequences of
    Theorem~\ref{Sigma}
    in view of
    Citation~\ref{Kriterium:Renz}:
    Let $\GenAbels$ be a normal subgroup $\BorelOver{\OkaSet}$ containing
    $\UniOver{\OkaSet}$. Such a subgroup always arises as a preimage
    of a subgroup $\GenSubTorus \leq\TorusOver{\OkaSet}$.
    A homomorphism in $\GrHom{\BorelOver{\OkaSet}}{\RRR}$ vanishes on
    $\GenAbels$ if and only if the corresponding linear form
    in $\KernDual$ vanishes on
    $\SubKernel:=\Tensor{\ZZZ}{\GenSubTorus}{\RRR} \leq
    \Kernel\isom\Tensor{\ZZZ}{\TorusOver{\OkaSet}}{\RRR}$.
    Thus, we are lead to consider
    $\SubDual := \SetOf{\linform\in\KernDual\suchthatvrule\linform(\SubKernel)=0}
    \leq\KernDual$. The codimension of $\SubDual$ in $\KernDual$ is nothing but
    the $\ZZZ$--rank of $\GenSubTorus$. According to
    Citation~\ref{Kriterium:Renz},
    $\GenAbels$ is of type \F{\type} provided the intersection
    $\BNSc{\type}{\BorelOver{\OkaSet}} \cap \SubDual$ is trivial.

    By the upper bound of Theorem~\ref{Sigma},
    $\BNSc{\type}{\BorelOver{\OkaSet}}$ is a subset of
    $\Conv{\type}{\NegativeRootsRestricted}$
    which is a finite union of $\type$--dimensional subspaces
    of $\KernDual$. Hence $\BNSc{\type}{\BorelOver{\OkaSet}}$
    has a trivial intersection with  $\SubDual$ provided
    the $\ZZZ$--rank of $\GenSubTorus$ is $\geq\type$ and 
    $\SubKernel$ (and hence $\SubDual$) is ``in general position''.

    The lower bound, analogously, implies that the $\ZZZ$--rank
    of $\GenSubTorus$ is at least $\type$ if
    $\GenAbels$ is of type \F{\type}. Hence in the ``generic case'',
    the finiteness length of $\GenAbels$ is given by the
    minimum of $\CardOf{\Set}-1$ and the $\ZZZ$--rank of
    $\GenSubTorus$.
  \end{rem}    
  \begin{rem}
    On the other side, it is not difficult to find subgroups of
    $\BorelOver{\OkaSet}$ containing $\UniOver{\OkaSet}$ that are
    not ``in general position''. Eg, consider the group
    $$
      \SetOf{
      \left(\begin{matrix}
        \RingEl & *            & *       \\
                & \RingEl^{-2} & *       \\
                &              & \RingEl
      \end{matrix}\right)
      \in\SlIndOver{3}{\OkaSet}
      \suchthatvrule
      \RingEl\in\UnitsOf{\OkaSet}}
      .
    $$
    In the number field case, a related group is discussed
    in~\cite{Abels.Brown:1987}.

    For $\CardOf{\Set}\geq2$, this group is finitely generated: we can
    write the elements in the upper right corner as products of
    commutators.  This implies that the upper bound of
    Theorem~\ref{Sigma} is not sharp.

    However,  independently of $\CardOf{\Set}$, this group is not
    finitely presented: It is step--$2$--nilpotent-by-abelian. If it was
    finitely presented then,
    by \cite[Corollary~5.8]{Bieri.Strebel:1980}, the subgroup in the
    upper right corner would be finitely generated as a normal subgroup.
    This is impossible since the torus acts trivially on it.
    It follows that the lower bound of
    Theorem~\ref{Sigma} is not sharp, as well.
  \end{rem}

  We now embark on the proof of Theorem~\ref{Sigma}. Both inclusions
  will be proved separately. The lower bound is derived in Lemma~\ref{UntereSigmaCompl}
  and the upper bound in Corollary~\ref{ObereSigmaCompl}.
  \begin{lemma}\label{UntereSigmaCompl}
    For $\type<\CardOf{\Set}$, we have the inclusion:
    \[
      \ConvSub{\type}{\Set}{\BaseRootsRestricted}
      \subseteq
      \BNSc{\type}{\BorelOver{\OkaSet}}
      .
    \]
  \end{lemma}
  \begin{proof}
    Let
    $
     \linform = 
       \Sum{\place\in\Set}{%
         \KoeffizientAt{\place}\AffRootAtInd{\place}{\IndexAt{\place}}
       }
    $
    be a linear form  on $\Kernel$. Suppose $\type$ coefficients
    are greater than $0$ whereas the others vanish.
    We will apply the technique of Section~\ref{GeometrischeVersion} to
    prove that the system
    $
      \FamOf{%
        \FundamentalGroup{\type-1}{%
          \PreIm{(\linform \circ \Restr{\ProProd}{\BorelSpace})}{%
            [\real,\infty)
          }
        }
      }{\real\in\RRR}
    $
    is not essentially trivial.

    In Section~\ref{GeometrischeVersion},
    we have seen how to construct a tree $\TreeAt{\place}$ 
    associated to a base root $\AffRootAtInd{\place}{\IndexAt{\place}}$
    that is a retract of the building $\AffBuildAt{\place}$ and
    that carries a height function
    $\TreeHightAt{\place} \mapcolon \TreeAt{\place} \rightarrow \RRR$
    related to the height function
    $
      \AffHightAt{\place}
      = \AffRootAtInd{\place}{\IndexAt{\place}} \circ \ProAt{\place} \mapcolon 
      \AffBuildAt{\place} \rightarrow \RRR
    $
    via the retraction map---see Lemma~\ref{BaumRetrakt} and
    Observation~\ref{SpezielleApartments}.

    Put
    $
      \RestrTreeProd :=
      \CrossProd{\KoeffizientAt{\place}\neq 0}{%
        \TreeAt{\place}
      }
    $.
    The induced retraction
    $
      \CrossProd{\KoeffizientAt{\place}\neq 0}{%
        \AffBuildAt{\place}
      }
      \rightarrow \RestrTreeProd
    $
    has a section
    $
      \RestrTreeSection \mapcolon
      \RestrTreeProd \rightarrow
      \CrossProd{\KoeffizientAt{\place}\neq 0}{%
        \AffBuildAt{\place}
      }
    $.
    Since $\type<\CardOf{\Set}$ there is at least one place
    $\place\in\Set$ with
    $\KoeffizientAt{\place}=0$. Hence we find a continuous map
    $
      \RestrTreeCorrection \mapcolon
      \RestrTreeProd \rightarrow
      \CrossProd{\KoeffizientAt{\place}=0}{%
        \AffAppAt{\place}
      }
    $
    such that, for each point $\RestrTreeProdPoint\in\RestrTreeProd$,
    the tuple
    $(\RestrTreeCorrectionEv{\RestrTreeProdPoint},\RestrTreeSectionEv{\RestrTreeProdPoint})
    \in 
    \CrossProd{\KoeffizientAt{\place}=0}{\AffAppAt{\place}}
    \times
    \CrossProd{\KoeffizientAt{\place}\neq 0}{\AffBuildAt{\place}}
    \subseteq
    \AffBuildProd$
    belongs to $\BorelSpace=\PreImOf{\ProProd}{\Kernel}$.
    Thus, we see that $\RestrTreeProd$ is a retract of $\BorelSpace$.

    Hence,
    $\Sum{\place\in\Set}{\KoeffizientAt{\place}\TreeHightAt{\place}}$
    is a height function on $\RestrTreeProd$ compatible with the height
    function induced by $\linform$ on $\BorelSpace$, and the following diagram
    commutes:
    \ifXYpic
    \[
      \xymatrix{
        {\BorelSpace}
        \ar@{->>}@<0.75mm>[r]
        \ar[d]_{\linform\circ\Restr{\ProProd}{\BorelSpace}}
      &
        {\RestrTreeProd}
        \ar@{ >->}@<0.75mm>[l]
        \ar[d]^{
          \Sum{\place\in\Set}{\KoeffizientAt{\place}\TreeHightAt{\place}}
        }
      \\
        {\RRR}
        \ar@{=}[r]
      &
        {\RRR}
      \\
      }
    \]
    \else
    \[
      \begin{array}{ccc}
        \BorelSpace & \retraktarrows & \RestrTreeProd\\
        \ldec{\linform\circ\Restr{\ProProd}{\BorelSpace}} 
        && 
        \ldec{\Sum{\place\in\Set}{\KoeffizientAt{\place}\TreeHightAt{\place}}}
        \\
        \RRR & = & \RRR
      \end{array}
    \]
    \fi
    Thus, arguing as in the proof of Lemma~\ref{Lemma:Non_Trivial}, we see that
    the system
    \(
      \FamOf{%
        \FundamentalGroup{\type-1}{%
          \PreIm{(\linform \circ \Restr{\ProProd}{\BorelSpace})}{%
            [\real,\infty)
          }
        }
      }{\real\in\RRR}
    \)
    of homotopy groups is not essentially trivial since the
    corresponding system
    \[
      \bbFamOf{%
        \bbRedHomTo{\type-1}{%
          \bbSetOf{(\TreePointAt{\place})\in\RestrTreeProd
            \bbsuchthatvrule
            \real \leq \Sum{\place\in\Set}{\KoeffizientAt{\place}\TreeHightAt{\place}(\TreePointAt{\place})}
          }
        }
      }{\real\in\RRR}
    \]
    of reduced homology groups
    on the retract $\RestrTreeProd$ is not essentially
    trivial by Lemma~\ref{BaumProdukt}.
  \end{proof}

  We turn to the upper bound.
  \begin{observation}\label{Einschraenkung}
    Let $\SigmaSet$ be an $\type$--tame set of linear forms on $\Kernel$
    and $\linform$ be a linear form on $\Kernel$. A positive combination of
    elements of $\SigmaSet$ vanishes on $\Kern{\linform}$ if and only if
    it belongs to the span of $\linform$ in the dual $\KernDual$.
    Hence, $\SigmaSet$ restricts
    to an $\type$--tame set of linear forms on $\Kern{\linform}$ if and
    only if $\linform\not\in\Conv{\type}{\SigmaSet}$ and
    $-\linform\not\in\Conv{\type}{\SigmaSet}$.\qed
  \end{observation}
  The following lemma, which is the key observation for the upper
  bound of Theorem~\ref{Sigma}, is a modification
  of~\cite[Lemma~7.3]{Bux:1997}.
  \begin{lemma}\label{Entzwei}
    Let $(\SheetCompl,\Eukl,\SheetProj)$ be a $\SheetSet$--complex
    and $\linform$ be a linear form on $\Eukl$. Suppose
    $\SheetSet$ is $\type$--tame and
    $\linform\not\in\Conv{\type}{\SheetSet}$. Then the preimage
    $\SheetPreim:=\PreIm{(\linform\circ\SheetProj)}{[0,\infty)}$
    is
    $(\type-1)$--connected.
  \end{lemma}
  \proof
    Let
    $\FamOf{\SheetInd{\Index}}{\Index \in \NNN}$ be an increasing
    sequence in $\SheetCompl$. We will show that for each
    $\LowInd\in\NNN$, the space
    $$
      \SheetFiltrInd{\LowInd} := \Union{\Index\leq\LowInd}
                            \SheetPreim \cap \SheetInd{\Index}
    $$
    is $(\type-1)$--connected. Since any compact subset of
    $\SheetCompl$ is covered by some increasing sequence, it follows
    that $\SheetPreim$ is $(\type-1)$--connected.

    Observe that the case $\LowInd=1$ is
    trivial: $\SheetFiltrInd{1}=\SheetInd{1}$ is contractible.
    For $\LowInd>1$, assume by induction that
    $\SheetFiltrInd{\LowInd-1}$ is $(\type-1)$--connected.
    We obtain $\SheetFiltrInd{\LowInd}$ from
    $\SheetFiltrInd{\LowInd-1}$ by gluing a closed convex set
    $\SheetNewPiece\subseteq\SheetInd{\LowInd}$
    onto $\SheetFiltrInd{\LowInd-1}$:
    Put $\AffHalfSpace:=\PreIm{\linform}[0,\infty)$ and
    $\EuKern := \Kern{\linform}=\Boundary{\AffHalfSpace}$. We identify
    the sheet $\SheetInd{\LowInd}$ with $\Eukl$ via the projection
    $\SheetProj$. This way,
    $\SheetNewPiece\cong\AffHalfSpace\cap\ConvexStep$ for a $\SheetSet$--convex set
    $\ConvexStep\subseteq\Eukl$. The ``new
    piece'' $\SheetNewPiece$ is glued in along the set
    $\AffHalfSpace\cap\Boundary{\ConvexStep}$.
    This process does not add material unless $\EuKern$ contains an interior
    point of $\ConvexStep$ which we shall henceforth assume.
    In this case,
    $\EuKern\cap\Boundary{\ConvexStep}=
    \Boundary{(\AffHalfSpace\cap\ConvexStep)}$.

    Using~\cite[Theorem~6.4]{Bux:1997} and its terminology, we
    distinguish two cases:

    \begin{case}
      \item[$\SheetNewPiece$ is decomposable]
        $\SheetNewPiece=\ConvexStep\cap\AffHalfSpace$ is a product
        of a convex compact set and a Euclidean space, which is also
        a factor of $\AffHalfSpace$.
        As we can split off this Euclidean factor, we may assume without loss
        of generality that
        $\SheetNewPiece$ is compact. Then, the intersection
        $\SheetNewPiece\cap\EuKern$ is a ``free face'' of $\SheetNewPiece$,
        which we can push into---recall that
        $\EuKern$ contains an interior point of $\ConvexStep$. Hence,
        $\SheetNewPiece$ can be collapsed onto
        $\Boundary{\ConvexStep}\cap\AffHalfSpace$, which we thus recognize as a
        retract. Since this is the part along which $\SheetNewPiece$ is
        glued in, $\SheetFiltrInd{\LowInd}$ is homotopy equivalent to
        $\SheetFiltrInd{\LowInd-1}$.
      \item[$\SheetNewPiece$ is retractable]
        In this case, gluing in $\SheetNewPiece$ amounts to glue in
        $\Boundary{\SheetNewPiece}$ along
        $\AffHalfSpace\cap\Boundary{\ConvexStep}$. However, we can retract everything
        toward $\EuKern$. Thus, we could
        equivalently glue in
        $\Boundary{\SheetNewPiece}\cap\EuKern=\ConvexStep\cap\EuKern$
        along $\Boundary{\ConvexStep}\cap\EuKern=
        \Boundary{(\ConvexStep\cap\EuKern)}$. So we have to determine the
        homotopy type of the pair
        \(
          (\ConvexStep\cap\EuKern,\Boundary{(\ConvexStep\cap\EuKern)})
          .
        \)

        Again, we have to distinguish two cases. Let
        $\SheetSubSet\subseteq\SheetSet$ be the set of all linear forms
        in $\SheetSet$ that correspond to hyperplanes supporting $\ConvexStep$.
        \begin{case}
          \item[$-\linform\not\in\Conv{\type}{\SheetSubSet}$]
            We have seen in Observation~\ref{Einschraenkung} that
            in this case,
            the system $\SheetSubSet$ induces an $\type$--tame system of
            linear forms on $\EuKern$. Hence, $\ConvexStep\cap\EuKern$ is
            either retractable or decomposable with a compact factor of
            dimension $\geq\type$. The bound on the dimension follows
            from~\cite[Corollary~6.7]{Bux:1997}. In either case, the
            gluing does not change homotopy groups in dimensions $<\type$,
            ie, we have
            \(
              \FundamentalGroup{\Index}{\SheetFiltrInd{\LowInd}}
              =
              \FundamentalGroup{\Index}{\SheetFiltrInd{\LowInd-1}}
            \)
            for $0\leq\Index<\type$.
          \item[$-\linform\in\Conv{\type}{\SheetSubSet}$]
            In this case, $\linform$ is bounded from above on $\ConvexStep$.
            Since the set $\SheetNewPiece$ is retractable, it contains an
            infinite ray without its antipode. This ray cannot depart
            from $\EuKern$ since otherwise $\linform$ was unbounded on
            $\SheetNewPiece\subseteq\ConvexStep$. Hence the ray is parallel to
            $\EuKern$. Thus, the set $\ConvexStep\cap\EuKern$ is
            retractable as a subset of $\EuKern$: it contains a
            parallel ray. Hence gluing in this set does not change the
            homotopy type of $\SheetFiltrInd{\LowInd-1}$.~\hfill\qedsymbol
        \end{case}
      \end{case}
  
    \begin{figure}[h]\small
      \centering
      \begin{tabular}{c@{\hspace{1cm}}c@{\hspace{1cm}}c}
                             \CaseOne &                                \CaseTwoAlpha &                             \CaseTwoBeta \\
                                      & $-\linform\not\in\Conv{\type}{\SheetSubSet}$ & $-\linform\in\Conv{\type}{\SheetSubSet}$ \\
        $\SheetNewPiece$ compact      &                \multicolumn{2}{c}{$\SheetNewPiece$ retractable}                         \\
      \end{tabular}
    \end{figure}

  \begin{cor}\label{ObereSigmaCompl}
    For $\type < \CardOf{\Set}$, we have the inclusion
    $\BNSc{\type}{\BorelOver{\OkaSet}}\subseteq
    \Conv{\type}{\NegativeRootsRestricted}$.
  \end{cor}
  \begin{proof}
    By Observation~\ref{BorelSpace}, $\TupelOf{\BorelSpace,\Kernel,\ProProd}$
    is a $\NegativeRootsRestricted$--complex. In the proof of Lemma~\ref{Urbild:Kernel:Zusammenhaengend},
    we saw that $\NegativeRootsRestricted$ is $\CardOf{\Set}$--tame. Thus the claim follows
    from
    Lemma~\ref{Entzwei} in view of Remark~\ref{Sufficient}.
  \end{proof}
  Together with Lemma~\ref{UntereSigmaCompl} this proves
  Theorem~\ref{Sigma}\qed\newpage

\begin{refferences}

\def\selectlanguage#1{\relax}

  \newcommand{\BielPre}[2]{%
    Preprint series #1--#2, SFB~343; Bielefeld (19#1)%
  }%

  \book{Abel87}{Abels:1987}%
    \au{H Abels}
    \ti{Finite presentability of $S$--arithmetic groups -- Compact
        presentability of solvable groups}
    \lo{Springer LNM~1261 (1987)}{}

  \article{AbBr87}{Abels.Brown:1987}
    \au{H Abels} \au{K\,S Brown}
    \ti{Finiteness Properties of Solvable $S$--Arithmetic Groups: An Example}
    \lo{Journal of Pure and Applied Algebra 44 (1987)}{77--83}

  \article{AbTi97}{Abels.Tiemeier:1997}%
    \au{H Abels} \au{A Tiemeier}
    \ti{Compactness Properties of Locally Compact Groups}
    \lo{Transformation Groups 2 (1997)}{119--135}
    (see also \cite{Tiemeier:1997})

  \article{\AA{}ber86}{Aberg:1986}%
    \au{H \AA{}berg}
    \ti{Bieri--Strebel Valuations (of Finite Rank)}
    \lo{Proceedings of the London Mathematical Society 52 (1986)}{269--304}

  \notes{Abra94}{Abramenko:1994}
    \au{P Abramenko}
    \ti{{\selectlanguage{german}Chevalleygruppen}}
    \lo{Manuscript (Frankfurt, WS~1994/95)}{}

  \book{Abra96}{Abramenko:1996}%
    \au{P Abramenko}
    \ti{Twin Buildings and Applications to $S$--Arithmetic Groups}
    \lo{Springer LNM~1641 (1996)}{}

  \article{BBN59}{Baumslag.Boone.Neumann:1959}%
    \au{G Baumslag} \au{W\,W Boone} \au{B\,H Neumann}
    \ti{Some Unsolvable Problems About Elements and Subgroups of Groups}
    \lo{Mathematica Scandinavica 7 (1959)}{191--201}

  \article{Behr69}{Behr:1969}%
    \au{H Behr}
    \ti{{\selectlanguage{german}Endliche Erzeugbarkeit arithmetischer
         Gruppen \"uber Funktionenk\"orpern}}
    \lo{Inventiones Mathematicae 7 (1969)}{1--32}

  \preprint{Behr92}{Behr:1992}%
    \au{H Behr}
    \ti{Arithmetic groups over function fields;
        A complete characterization of finitely generated and
        finitely presented arithmetic subgroups of reductive algebraic
        groups}
     \lo{\BielPre{92}{033}}{};
     part of this is published in \cite{Behr:1998}

  \article{Behr98}{Behr:1998}%
    \au{H Behr}
    \ti{Arithmetic groups over function fields I;
        A complete characterization of finitely generated and
        finitely presented arithmetic subgroups of reductive algebraic
        groups}
     \lo{{\selectlanguage{german}(Crelles) Journal f\"ur die reine und angewandte Mathematik} 495 (1998)}{79--118}

  \article{BeBr97}{Bestvina.Brady:1997}%
    \au{M Bestvina} \au{N Brady}
    \ti{Morse theory and finiteness properties of groups}
    \lo{Inventiones mathematicae 129 (1997)}{445--470}

  \article{Bier99}{Bieri:1999}%
    \au{R Bieri}
    \ti{Finiteness length and connectivity length for groups}
    from: ``Geometric Group Theory down under'', (J Cossey, C\,F Miller,
    W Neumann, M Shapiro, editors),
    Proceedings of the special year in geometric group theory,
    Australia National University, Canbera, 1996, deGruyter (1999) 9--22 

  \article{BiEc74}{Bieri.Eckmann:1974}
    \au{R Bieri} \au{B Eckmann}
    \ti{Finiteness Properties of Duality Groups}
    \lo{Commentarii Mathematici Helvetici 49 (1974)}{74--83}

  \preprint{BiGe02}{Bieri.Geoghegan:2002}%
    \au{R Bieri} \au{R Geoghegan}
    \ti{Connectivity properties of group actions on non-positively
        curved spaces}
    \lo{Memoirs of the American Mathematical Society 161 (2003) no.\,765}{}

  \article{BiRe88}{Bieri.Renz:1988}%
    \au{R Bieri} \au{B Renz}
    \ti{Valuations on free resolutions and higher geometric
        invariants of groups}
    \lo{Commentarii Mathematicae Helvetici 63 (1988)}{464--497}

  \article{BiSt80}{Bieri.Strebel:1980}
    \au{R Bieri} \au{R Strebel}
    \ti{Valuations and Finitely Presented Metabelian Groups}
    \lo{Proceedings of the London Mathematical Society 41 (1980)}{439--464}

  \book{Bore91}{Borel:1991}%
    \au{A Borel}
    \ti{Linear Algebraic Groups}
    \edition{Second Enlarged Edition}
    \lo{Springer GTM 126, New York (1991)}{}

  \article{BoSe76}{Borel.Serre:1976}%
    \au{A Borel} \au{J-P Serre}
    \ti{{\selectlanguage{french}Cohomologie d'immeubles et de groupes S-arithm\'{e}tiques}}
    \lo{Topology 15 (1976)}{211--232}

  \book{Brow82}{Brown:1982}%
    \au{K\,S Brown}
    \ti{Cohomology of Groups}
    \lo{Springer GTM 87, New York (1982)}{}

  \article{Brow87}{Brown:1987}%
    \au{K\,S Brown}
    \ti{Finiteness Properties of Groups}
    \lo{Journal of Pure and Applied Algebra 44 (1987)}{45--75}

  \book{Brow89}{Brown:1989}
    \au{K\,S Brown}
    \ti{Buildings}
    \lo{Springer, New York, Berlin (1989)}{}

  \report{BrTi72}{Bruhat.Tits:1972}%
    \au{F Bruhat} \au{J Tits}
    \ti{{\selectlanguage{french}Groupes R\'{e}ductives sur un Corps Local I}}
    \lo{{\selectlanguage{french}Publications Math\'{e}matiques IHES 41 (1972)}}{5--252}

  \report{BrTi84}{Bruhat.Tits:1984}%
    \au{F Bruhat} \au{J Tits}
    \ti{{\selectlanguage{french}Groupes R\'{e}ductives sur un Corps Local II}}
    \lo{{\selectlanguage{french}Publications Math\'{e}matiques IHES 60 (1984)}}{5--184}

  \article{Bux97a}{Bux:1997}%
    \au{K-U Bux}
    \ti{Finiteness Properties of some Metabelian $S$--Arithmetic Groups}
    \lo{Proceedings of the London Mathematical Society 75 (1997)}{308--322}

  \thesis{Bux97b}{Bux:1997b}%
    \au{K-U Bux}
    \ti{{\selectlanguage{german}Endlichkeitseigenschaften aufl\"osbarer
        arithmetischer Gruppen \"uber Funktionenk\"orpern}}
    \lo{PhD thesis (Frankfurt, 1997)}{}

  \article{Bux02}{Bux:2002}%
    \au{K-U Bux}
    \ti{Finiteness Properties Soluble $S$--Arithmetic Groups -- A Survey}
    \lo{To appear: Conference Proceedings ``Groups: Geometric and Combinatorial Aspects'' Bielefeld 1999}{}

  \book{CaFr67}{Cassels.Froehlich:1967}
    \ed{J\,W\,S Cassels} \au{A Fr\"ohlich}
    \ti{Algebraic Number Theory. Proceedings of an Instructional
        Conference organized by the London Mathematical Society
        (A Nato Advanced Study Institute)
        with Support of the International Mathematical Union}
    \lo{Academic Press, London (1967)}{}

  \article{Chev60}{Chevalley:1960}%
    \au{C Chevalley}
    \ti{{\selectlanguage{french}Certain Sch\'{e}mas de Groupes Semi-Simples}}
    \lo{{\selectlanguage{french}S\'{e}minaire Bourbaki 13e ann\'{e}e (1960/61) n$\circ$~219}}{1--16}

  \book{GrRo84}{Gruenberg.Roseblade:1984}
    \ed{K\,W Gruenberg} \au{J\,E Roseblade}
    \ti{Group Theory: Essays for Philip Hall}
    \lo{Academic Press (1984)}{}

  \thesis{Koch97}{Kochloukova:1997}
    \au{D Kochloukova}
    \ti{The FP${}^m$--Conjecture for a class of metabelian groups}
    \lo{PhD thesis (Cambridge, 1997)}{}

  \article{Mein96}{Meinert:1996}%
    \au{H Meinert}
    \ti{The homological invariants for metabelian groups of
        finite Pr\"ufer rank: a proof of the $\Sigma^m$--conjecture}
    \lo{Proceedings of the London Mathematical Society 72 (1996)}{385--424}

  \article{Mein97}{Meinert:1997}%
    \au{H Meinert}
    \ti{Actions on $2$--complexes and the homotopical invariant
        $\Sigma^2$ of a group}
    \lo{Journal of Pure and Applied Algebra 119 (1997)}{297--317}

  \article{Neum37}{Neumann:1937}%
    \au{B\,H Neumann}
    \ti{Some remarks on infinite groups}
    \lo{Journal of the London Mathematical Society 12 (1937)}{120--127}
  \thesis{Renz88}{Renz:1988}
    \au{B Renz}
    \ti{{\selectlanguage{german}Geometrische Invarianten und Endlichkeitseigenschaften von
        Gruppen}}
    \lo{PhD thesis (Frankfurt, 1988)}{}

  \book{Rona89}{Ronan:1989}%
    \au{M Ronan}
    \ti{Lectures on Buildings}
    \lo{Perspectives in Mathematics,
        Academic Press, San~Diego, London (1989)}{}

  \book{Serr80}{Serre:1980}%
    \au{J-P Serre}
    \ti{Trees}
    \lo{Springer, Berlin Heidelberg New York (1980)}{}

  \notes{Stei67}{Steinberg:1967}%
    \au{R Steinberg}
    \ti{Lectures on Chevalley Groups}
    \lo{Mimeographed Notes, Yale University (1967)}{}

  \article{Stre84}{Strebel:1984}%
    \au{R Strebel}
    \ti{Finitely Presented Solvable Groups}
    \lo{from: \cite{Gruenberg.Roseblade:1984}}{257--314}

  \article{Stuh76}{Stuhler:1976}%
    \au{U Stuhler}
    \ti{{\selectlanguage{german}Zur Frage der endlichen Pr\"asentiertheit gewisser
        arithmetischer Gruppen im Funktionenk\"orperfall}}
    \lo{{\selectlanguage{german}Mathematische Annalen 224 (1976)}}{217--232}

  \article{Stuh80}{Stuhler:1980}%
    \au{U Stuhler}
    \ti{Homological properties of certain arithmetic groups in
        the function field case}
    \lo{Inventiones mathatematicae 57 (1980)}{263--281}

  \article{Tiem97}{Tiemeier:1997}%
    \au{A Tiemeier}
    \ti{A Local-Global Principle for Finiteness Properties of
      $S$--Artihmetic Groups over Number Fields}
    \lo{Transformation Groups 2 (1997)}{215--223}
    (see also \cite{Abels.Tiemeier:1997})

  \article{Tits87}{Tits:1987}%
    \au{J Tits}
    \ti{Uniqueness and Presentations of Kac-Moody Groups over Fields}
    \lo{Journal of Algebra~105 (1987)}{542--573}

  \article{Wall65}{Wall:1965}
    \au{C\,T\,C Wall}
    \ti{Finiteness Conditions for CW-Complexes}
    \lo{Annals of Mathematics 81 (1965)}{56--69}

  \book{Wall79}{Wall:1979}
    \ed{C\,T\,C Wall}
    \ti{Homological Group Theory}
    \lo{LMS Lecture Notes 35, Cambridge University Press (Cambridge 1979)}{}

  \book{Weil73}{Weil:1973}%
    \au{A Weil}
    \ti{Basic Number Theory}
    \edition{Reprint of the Second Edition (1973)}
    \lo{Springer Classics in Mathematics, Berlin Heidelberg (1995)}{}

\end{refferences}

\end{document}